%% file: main.tex
\documentclass[a4,11pt]{paper}
\usepackage{enumitem,color,amssymb,bm,amsmath,tikz-cd,mathtools}
\usepackage[english]{babel}  
\usepackage{graphicx}
\usepackage{lipsum}
\usepackage{xcolor}
\usepackage{booktabs}
\usepackage{scalerel}
\usepackage{ntheorem}
{\theoremstyle{empty}
	}
\usepackage[top=1in, left=1in, right=1in, bottom=0.8in]{geometry}
\newcommand{\R}{\mathbb{R}}
\newcommand{\Z}{\mathbb{Z}}

\newcommand{\Abb}{{\mathbb A}}
\newcommand{\Ical}{{\mathcal I}}
\newcommand{\Jcal}{{\mathcal J}}
\newcommand{\spinct}{{\mathfrak t}}
\newcommand{\Dcal}{{\mathcal D}}
\newcommand{\Pcal}{{\mathcal P}}

\newcommand{\Tri}{\mathsf{Tri}}
\newcommand{\tdelta}{\widetilde{\delta}}

\newcommand{\z}{\boldsymbol{z}}
\newcommand{\w}{\boldsymbol{w}}
\newcommand{\x}{\mathbf{x}}

\newcommand{\zs}{{\vec{z}}}
\newcommand{\ws}{{\vec{w}}}

\newcommand{\alphas}{{\boldsymbol{\alpha}}}
\newcommand{\betas}{{\boldsymbol{\beta}}}
\newcommand{\gammas}{\boldsymbol{\gamma}}

\newcommand{\SpinC}{{\mathrm{Spin}^c}}

\newcommand{\p}{\mathbf{p}}
\def\endproof{\relax\ifmmode\expandafter\endproofmath\else
	\unskip\nobreak\hfil\penalty50\hskip.75em\hbox{}\nobreak\hfil\bull
	{\parfillskip=0pt \finalhyphendemerits=0 \bigbreak}\fi}
\def\endproofmath$${\eqno\bull$$\bigbreak}
\def\bull{\vbox{\hrule\hbox{\vrule\kern3pt\vbox{\kern6pt}\kern3pt
			\vrule}\hrule}}

\newtheorem{thm}{\bfseries\color{black!50!blue} Theorem}[section]
\newtheorem{prop}[thm]{\bfseries\color{black!50!blue} Proposition}
\newtheorem{cor}[thm]{\bfseries\color{black!50!blue} Corollary}

\newtheorem{lem}[thm]{\bfseries\color{black!50!blue} Lemma}

\newtheorem{rem}[thm]{\bfseries\color{black!50!blue} Remark}

\newtheorem{defi}[thm]{\bfseries\color{black!50!blue} Definition}

\newcommand{\cptwo}{{$\CPbar$}}
\newcommand{\spinc}{{\mathfrak s}}
\newcommand{\sps}{\mathfrak{s}}
\newcommand{\spt}{\mathfrak{t}}
\newcommand{\cfm}[1]{\mathrm{CFL}^-_{#1}}
\newcommand{\cfp}[1]{\mathrm{CFL}^\circ_{#1}}

\newcommand{\hfm}[1]{\mathrm{HFL}^-_{#1}}
\newcommand{\hfp}[1]{\mathrm{HFL}^\circ_{#1}}
\newcommand{\hfs}[1]{\mathrm{HFL}^\star_{#1}}

\newcommand{\ftwo}{\mathbb{F}}
\newcommand{\uu}{\mathrm{u}}
\newcommand{\vv}{\mathrm{v}}

\newcommand{\gen}{\boldsymbol{\mathrm{S}}}
\newcommand{\gene}[1]{\boldsymbol{\mathrm{#1}}}

\newcommand{\chai}{\mathcal{C}}
\newcommand{\link}{\mathcal{L}}

\newcommand{\Goo}{Admissible }

\newcommand{\Nic}{Nice }

\newcommand{\heeg}{\mathcal{H}}
\newcommand{\hee}{H}
\newcommand{\doma}{\mathcal{D}}
\newcommand{\tria}{\mathcal{T}}

\subsectionfont{\large\sf\bfseries\color{black!50!blue}} 
\sectionfont{\Large\sf\bfseries\color{black!50!blue}} 
\titlefont{\LARGE\sf\bfseries\color{black!50!blue}} 
\numberwithin{equation}{section}
\newif\ifshowrevisions
\showrevisionstrue

\newcommand{\CPbar}{\overline{\mathbb{CP}}^{2}}
\newcommand{\grw}{\operatorname{gr}_{\w}}
\newcommand{\grz}{\operatorname{gr}_{\z}}
\newcommand{\PD}{\operatorname{PD}}

\begin{document}
	\title{Grid triples, link cobordisms in \cptwo \ and the computation of Heegaard-Floer TQFT}
	\author{Roohallah Mahkam}
	\institution{\sf{School of Mathematics, Institute for Research in Fundamental Sciences (IPM), Tehran, Iran}}
	%————————————————————————————————
	\maketitle

\begin{abstract}
	We give a combinatorial description of the Heegaard--Floer
	link-cobordism maps associated with pointed triple grid diagrams in
   the complement of three $4$-balls in complex projective plane with reverse orientation.
	
	Every contributing index-zero holomorphic triangle decomposes into
	Euclidean triangular sheets. We characterize the compatible collections
	of sheets that contribute by an explicit geometric separation criterion,
	expressed as a lattice-point test, and prove that each has mod-\(2\)
	count one. Individual sheets may wrap around the torus in the unblocked
	theory, whereas they are embedded in the blocked theory. In each fixed
	\(\SpinC\) structure, the map is determined by a finite enumeration.
	
	A difference number computed from the triangular regions determines
	the \(\SpinC\) structure of a Whitney triangle. Together with the
	topology of the associated surface, this gives an explicit formula
	for the absolute bigrading of each summand of the cobordism map.
	After basepoint reduction, the formula determines the Laurent monomial
	representing any nonzero map on the rank-one torsion-free two-variable
	knot Floer modules, yielding vanishing and ideal-containment criteria.
	
	We apply the construction to an affine family of grid links
	\(K(n;p,s)\), obtaining explicit cobordisms involving infinite families
	of torus knots and links, together with closed formulas for the
	bigradings of their distinguished grid states.
\end{abstract}

\section{Introduction}
A recurring theme in Heegaard--Floer theory is the passage from invariants of individual three-manifolds or knots to maps associated with cobordisms.  Ozsv\'ath and Szab\'o's original construction associated Floer groups to closed three-manifolds and, through knot Floer homology, to knots in three-manifolds \cite{OS-3m1, OS-3m2, OS-knot}; Rasmussen's early work provided a closely related viewpoint on the knot filtration \cite{Ras} and Eftekhary extended the definition to knots that are not null-homologous \cite{Ef-LFH}.  
The theory rapidly developed into a collection of invariants sensitive to surfaces, concordance, the Thurston norm, and four-dimensional topology; see, for example, \cite{OS-tau,Ni-Thurston-norm}.  From this perspective, understanding the maps induced by cobordisms is not secondary to the theory: it is part of the mechanism by which the invariants interact with topology.

A second fundamental development concerned naturality.  Heegaard Floer groups are originally defined from auxiliary Heegaard diagrams, and constructing a TQFT requires canonical control over the identifications associated with changes of diagram.  Juh\'asz established functoriality for sutured cobordisms \cite{Juhasz-naturality}, while the naturality problem for Heegaard Floer homology and the associated mapping class group actions was developed systematically by Juh\'asz, Thurston and Zemke \cite{JTZ} and independently by Alishahi and Eftekhary \cite{AE-2}.  More recently, a decoration-independent version of the theory was constructed in \cite{Ef-TQFT}, producing a Heegaard--Floer TQFT for oriented link cobordisms in which the auxiliary decoration is no longer part of the final invariant. These results made it possible to regard the Floer groups as genuinely functorial objects rather than groups defined only up to unspecified isomorphism. These theories have since been used in a variety of four-dimensional applications, including concordance maps, slice disks, ribbon concordance, bridge-type questions and different notions of unknotting bounds \cite{JM-Concordance,JZ-slice-disks,Zemke-Ribbon,JMZ, AE-unknotting, Ef-RTR, Ef-torsion}.

Geometric Floer constructions, which are by nature deeply analytic,  become truly computational only after the analytic data have been translated into sufficiently rigid combinatorial structures.
Grid diagrams provide perhaps the most successful example of such a translation.  Their use in knot theory goes back well before Floer homology; see, for instance, Cromwell's work on arc presentations \cite{Cromwell}.  In Heegaard Floer theory they became the basis of the combinatorial descriptions of knot and link Floer homology developed by Manolescu, Ozsv\'ath, Sarkar, Szab\'o and Thurston \cite{MOS,MOST}, and later organized systematically in \cite{OSS-book}.  In a grid diagram the generators become permutations, the relevant holomorphic disks become empty rectangles, and Maslov and Alexander gradings become finite counting functions.  This replacement of analysis by planar combinatorics was strong enough not only to compute the homology but also to support further constructions involving skein triangles, cobordism maps and related operations; see, among others, \cite{LMW,Ciprian-unoriented,Wong,Wong-2017}.  Sarkar's polygon index formulas \cite{Sarkar-index} provide another bridge between the analytic and combinatorial sides of the theory, and Eftekhary's comparison between sign assignments extends the equivalence to the theory with integer coefficients \cite{Ef-orientation}.

The question motivating this paper is whether a comparable phenomenon occurs for link-cobordism maps arising from Heegaard triples.  Ordinary grids are adapted to pairs of attaching-curve systems and therefore to differentials.  Triangle maps require three systems of curves simultaneously.  Triple grid diagrams supply precisely such a configuration.
Blackwell, Gay and Lambert-Cole used them to construct Lagrangian surfaces in \(\mathbb{CP}^2\) \cite{BGL-triple-grid}; see also Blackwell's thesis \cite{Bla22}.  Related lattice and hexagonal descriptions of curves and surfaces in \(\mathbb{CP}^2\) appear, for example, in \cite{Zupan-hexagonal}.  A triple grid consists of three families of curves of slopes \(0\), \(\infty\) and \(-1\) on a torus.  Any two families form an ordinary grid diagram, while the three families determine the oppositely oriented projective-plane cobordism
\[
\CPbar\setminus(B^4\sqcup B^4\sqcup B^4).
\]
Thus the same elementary diagram simultaneously records three boundary links, an embedded surface joining them, and the Heegaard triple from which the Floer cobordism map is defined.
This paper studies the Heegaard--Floer meaning of that diagram.

The basic regions of a triple grid are triangles $L_{i,j}$ and $U_{i,j}$ of two orientation types, as well as hexagons $H_{i,j}$ for $i,j=1,\ldots,n$.  A \emph{nice} domain is a compatible sum of \(n\) Euclidean triangular sheets satisfying a geometric separation condition: their graphs over a common source triangle are pairwise disjoint.   Each embedded triangle may of course be constructed as a union of basic regions. We give a completely combinatorial computation of the triangle map associated with the two-variable link Floer complex and its blocked specializations.

\begin{thm}\label{thm:intro-1}
	The holomorphic triangle map associated with a pointed triple grid of 
	size \(n\) is combinatorial in the split model: every contributing domain is nice,
	and every nice class has mod-\(2\) holomorphic
	count one. Only finitely many classes contribute in each fixed
		\(\SpinC\) summand. The corresponding surface has Euler
	characteristic \(-n\), and its genus is determined by its bridge graph and
	boundary component counts. 
\end{thm}

The first statement follows from the Maslov-index and fat-diagonal formula of Sarkar \cite{Sarkar-index}: we show that an index-zero contributing domain has zero intersection with the fat diagonal, so a holomorphic representative separates into \(n\) genuine triangular sheets.  
The separation condition is itself combinatorial. If the lifted sheets
are \(A_i+\lambda_iT\), where \(T\) is the unit right triangle and
\(\lambda_i\) is the signed leg length, it is equivalent to
\[
\bigl(A_i-A_j+(\lambda_i-\lambda_j)T\bigr)\cap(n\mathbb Z)^2
=\varnothing\qquad(i<j).
\]
Proposition~\ref{prop:sheet-separation} also computes the Maslov index by
counting lattice points in these difference triangles.

	The \(\SpinC\) structure associated with a triangle class can also be read directly
	from its domain. If the formal domain associated with a triangle class $\psi$ is
	\[
	\mathcal D(\psi)=\sum_{i,j}
	\bigl(u_{i,j}U_{i,j}+h_{i,j}H_{i,j}+l_{i,j}L_{i,j}\bigr),
	\]
	we define
	\[
	\delta(\psi)=\sum_{i,j}(l_{i,j}-u_{i,j})\qquad\text{and}\qquad
	\tdelta(\psi)=-\frac{\delta(\psi)}{2n}.
	\]
\begin{prop}\label{prop:intro-1}
		Let \(\heeg\) be a centrally marked appropriate triple grid of size \(n\),
		and let \(H\) be the generator of \(H_2(X_\heeg;\mathbb Z)=H_2(\CPbar;\Z)\). For every Whitney triangle class
		\(\psi\),
		\[
		\tdelta(\psi)\in\tfrac12+\mathbb Z,\qquad
		\bigl\langle c_1(\spinc(\psi)),H\bigr\rangle
		=2\tdelta(\psi)=-\frac{\delta(\psi)}n.
		\]
\end{prop}

Combining this topological information with the absolute grading formula for link-cobordism maps \cite{Zemke-gradings-revision} gives an explicit bidegree for every \(\SpinC\)-summand. 

\begin{prop}\label{prop:intro-2}
		For a centrally marked appropriate triple grid of size \(n\), the associated
		surface is relatively nullhomologous. With \(H^2=-1\) as before, and
		\(c_{\spinct}=\langle c_1(\spinct),H\rangle\), the unreduced
		two-variable map has cobordism bidegree
		\[
		\deg_{\w}(F_{\heeg,\spinct})
		=\deg_{\z}(F_{\heeg,\spinct})
		=\frac{3-2n-c_{\spinct}^{2}}4.
		\]
		For a triangle representing \(\spinct\),
		\(c_{\spinct}=2\tdelta(\psi)=-\delta(\psi)/n\), as in
		Proposition~\ref{prop:intro-1}.
\end{prop}

	For a knot \(K\subset S^3\), let \(\Ical(K)\) denote the torsion-free
	quotient of its doubly based two-variable knot Floer homology over
	\(\Abb=\mathbb F[\uu,\vv]\). Choose a monomial-ideal realization
	\(\Ical(K)\cong\Jcal(K)[\eta(K)]\), where the shift \(\eta(K)\)
	is added to the usual monomial bidegree of the monomial ideal $\Jcal(K)\subset \mathbb A$. The grading computation gives
	the following consequence of Theorem~\ref{thm:ideal-inclusion}.

\begin{thm}\label{thm:intro-ideals}
		Suppose that the three boundary links of \(\heeg\) are knots
		\(K_0,K_1,K_2\). After the \(\z\)-basepoint-action reduction, each
		\(\SpinC\) structure \(\spinct\) determines a decoration-independent map
		\[
		\overline F_{\heeg,\spinct}:\Ical(K_0)\otimes_{\Abb}\Ical(K_1)
		\longrightarrow\Ical(K_2).
		\]
		If this map is nonzero, its expression under the ideal realizations is
		multiplication by a Laurent monomial \(\uu^a\vv^b\), and
		\[
		\uu^a\vv^b\Jcal(K_0)\Jcal(K_1)\subset\Jcal(K_2).
		\]
		The exponents are determined by
		\[
		(-2a,-2b)
		=\left(\frac{1-c_{\spinct}^2}{4},
		\frac{1-c_{\spinct}^2}{4}-n+1\right)
		+\eta(K_0)+\eta(K_1)-\eta(K_2).
		\]
\end{thm}

We also study a concrete arithmetic family.  For integers \(n,p,s\) we define affine grid links \(K(n;p,s)\) by placing the two markings in each column according to affine permutations.  Their elementary symmetries, component numbers and distinguished grid gradings can be computed explicitly.  Under suitable invertibility assumptions the three boundary diagrams of the associated triple grid are again affine and take the form
\[
K(n;p,s),\qquad
K(n;1-p,-s)\qquad\text{and}\qquad
K(n;-pr,-sr),
\]
respectively, where \(r=(p-1)^{-1}\pmod n\).
The subfamily \(K(n;p,pq)=T_{q,n-pq}\)
recovers torus links and produces explicit cobordisms between families of torus knots and links.  These examples provide a testing ground in which the surface topology, the combinatorial triangle enumeration, the \(\SpinC\) labels and the grading constraints can all be calculated simultaneously.

The paper is organized as follows.  We first recall the multi-based link Floer complexes, decorated link cobordisms and their triangle maps.  We then introduce pointed triple grids and compute the topology of the associated surfaces.  The next sections develop nice domains and the difference-number description of \(\SpinC\) structures, and prove the combinatorial formula for the triangle maps, including the geometric separation test and finiteness in each fixed $\SpinC$ summand.  We then derive the absolute grading formula and its consequences for nonzero torsion-free summands.  Finally, we study affine grids and the resulting torus-link families.

{\color{black!50!blue}\textbf{Acknowledgements.}}
I am deeply grateful to my PhD advisor, Professor Eaman Eftekhary, for his kind and generous support throughout my doctoral research. His insight and guidance were essential to the results of this work, and his invaluable help in writing and editing the manuscript greatly improved its presentation.

\section{Heegaard-Floer TQFT for link cobordisms}

\subsection{Pointed links, link cobordisms and trisection diagrams}
Let $L$ be an oriented link in a closed $3$-manifold $Y$ together with a finite collection of points on it, with at least one point on every connected component of $L$. Considering $L$ as a $1$-manifold, we can choose neighborhoods around each point on $L$. These neighborhoods give us a collection of disjoint oriented arcs on $L$ with the orientations induced by $L$. Let $w_i$ and $z_i$ for $i= 1, 2,\dots, k$ be respectively the initial and terminal points of these arcs. By a pointed link $\link$ we mean a triple
$$\link = (L, \boldsymbol{w} = \{ w_1, \dots, w_k \}, \boldsymbol{z} = \{ z_1, \dots, z_k \})$$ defined in this way. To every pointed link $\link \subset Y$ one can associate a pointed Heegaard diagram $\hee = (\Sigma, \boldsymbol{\alpha},  \boldsymbol{\beta}, \boldsymbol{w}, \boldsymbol{z})$ 
where $\Sigma$ is a separating closed surface in $Y$, and $\boldsymbol{\alpha}, \boldsymbol{\beta}$ are collections of disjoint simple closed curves bounding disks on the two sides of $\Sigma$ which, together with $\Sigma$, cut $Y$ into disjoint 3-balls and 
$L \cap \Sigma = \boldsymbol{w} \cup \boldsymbol{z}$ where each of $\boldsymbol{w}$ and $\boldsymbol{z}$ contains precisely one basepoint in each connected component of $\Sigma - \boldsymbol{\alpha}$ and $\Sigma - \boldsymbol{\beta}$. Also, given a pointed Heegaard diagram $\hee = (\Sigma, \boldsymbol{\alpha},  \boldsymbol{\beta}, \boldsymbol{w}, \boldsymbol{z})$ with $\#\boldsymbol{w} = \#\boldsymbol{z}$, one can construct a sutured manifold that determines an oriented pointed link $\link_{\hee}$ in a closed $3$-manifold $Y_{\hee}$. For the multi-based Heegaard-diagram convention, see \cite[Definition~3.1]{Zemke-1}; for the sutured-manifold viewpoint, see \cite[Section~2.1]{AE-2} and \cite{AE-1}.

\begin{defi} 
	A \emph{pointed trisection diagram} consists of a Heegaard triple $\heeg = (\Sigma, \boldsymbol{\alpha}, \boldsymbol{\beta}, \boldsymbol{\gamma}, \bf{w}, \bf{z})$ where $\Sigma$ is a closed surface of genus $g$, and  $\boldsymbol{\alpha}$, $\boldsymbol{\beta}$ and $\boldsymbol{\gamma}$ are all collections of disjoint simple closed curves on $\Sigma$ with $|\boldsymbol{\alpha}| = |\boldsymbol{\beta}| = |\boldsymbol{\gamma}| \geq g$, and $(\boldsymbol{w}, \boldsymbol{z}) \subset \Sigma - \boldsymbol{\alpha} - \boldsymbol{\beta} - \boldsymbol{\gamma}$ are two sets of basepoints. Also, these collections of closed curves and basepoints are chosen so that $\Sigma - \boldsymbol{\alpha}$,  $\Sigma - \boldsymbol{\beta}$ and $\Sigma - \boldsymbol{\gamma}$ are genus-$0$ surfaces and have the same number of components so that each one of these components includes precisely one basepoint from $\boldsymbol{w}$ and one from $\boldsymbol{z}$.
\end{defi}

Given a pointed trisection diagram $\heeg = (\Sigma, \boldsymbol{\alpha}, \boldsymbol{\beta}, \boldsymbol{\gamma}, \bf{w}, \bf{z})$, define
$$\hee_0 = (\Sigma, \boldsymbol{\alpha},  \boldsymbol{\beta}, \boldsymbol{w}, \boldsymbol{z}) , \hee_1 = (\Sigma, \boldsymbol{\beta},  \boldsymbol{\gamma}, \boldsymbol{w}, \boldsymbol{z}) \ \mathrm{and}\  \hee_2 = (\Sigma, \boldsymbol{\alpha},  \boldsymbol{\gamma}, \boldsymbol{w}, \boldsymbol{z}),$$ then each $\hee_i$ is a Heegaard diagram for a pointed link $\link_i$ in a $3$-manifold $Y_i$. In the rest of this paper, we sometimes refer to pointed trisection diagrams simply as trisection diagrams.

Given a trisection diagram $\heeg$ one can construct a smooth cobordism $X_\heeg$ which is {a compact oriented $4$-manifold with boundary $-Y_0\amalg -Y_1\amalg Y_2$} and a smooth surface $S_\heeg \subset X_\heeg$ {with boundary $-\link_0\amalg -\link_1\amalg\link_2$} in the following way. Let $U_0$, $U_1$ and $U_2$ denote the compression bodies determined by the surface $\Sigma$ and $\boldsymbol{\gamma}$, $\boldsymbol{\alpha}$ and $\boldsymbol{\beta}$ respectively. The basepoints $\boldsymbol{w}$ and $\boldsymbol{z}$ determine a union of oriented arcs in $U_i$, which is denoted by $T_i$. We denote the end-points of $T_i$ on $\partial U_i = \Sigma$ by $\boldsymbol{w}_i$ and $\boldsymbol{z}_i$. Let $\triangle$ denote a triangle with vertices $v_0$, $v_1$, $v_2$ and edges $e_0$, $e_1$, $e_2$ so that $e_i$ connects $v_{i+1}$ to $v_{i-1}$ for every $i \in \Z /3 $. Define
$$X_\heeg = \frac{(\Sigma \times \triangle) \coprod (\coprod_{i \in \Z /3} U_i \times e_i)}{\Sigma \times e_i \sim \partial U_i \times e_i},  \ S_\heeg = \frac{((\boldsymbol{w} \cup \boldsymbol{z}) \times \triangle) \coprod (\coprod_{i \in \Z /3} T_i \times e_i)}{(\boldsymbol{w} \cup \boldsymbol{z}) \times e_i \sim (\boldsymbol{w}_i \cup \boldsymbol{z}_i) \times e_i}.$$

We use the standard Heegaard-triple orientation
\(\partial X_{\alpha\beta\gamma}=-Y_{\alpha\beta}\sqcup
-Y_{\beta\gamma}\sqcup Y_{\alpha\gamma}\):
over the interior of the parameter triangle its sheets are
\(\z\times\triangle\) and \(-\w\times\triangle\).

Then $X_\heeg$ will be a cobordism from $Y_0 \amalg Y_1$ to $Y_2$ and $S_\heeg$ will be a cobordism from $\link_0 \amalg \link_1$ to $\link_2$. We denote these cobordisms by 
$$(X_\heeg, S_\heeg): (Y_0 \amalg Y_1, \link_0 \amalg \link_1) \rightarrow (Y_2, \link_2).$$
Following \cite[Definitions~2.1--2.4]{Zemke-1}, a decoration of a link cobordism $(X,S)$ with multi-based boundary links consists of a properly embedded $1$-manifold $\mathcal A\subset S$ (the \emph{dividing set}) and a partition of the regions of $S\setminus\mathcal A$ into $S_{\w}$ and $S_{\z}$, which meet along $\mathcal A$. On each boundary link, the points of $\mathcal A\cap\partial S$ lie between consecutive basepoints, so that every component of $L_i\setminus\mathcal A$ contains exactly one basepoint; the $w$-basepoints lie on $S_{\w}$ and the $z$-basepoints on $S_{\z}$. Zemke also specifies a compatible coloring of the basepoints and the components of $S\setminus\mathcal A$. Here we use the two-color specialization that identifies all $w$-variables with $\uu$ and all $z$-variables with $\vv$. We write $\mathcal F=(S,\mathcal A)$ for this decorated surface, with its coloring understood.

The formula above for $S_\heeg$ specifies the underlying surface and its boundary basepoints. In order to specify the dividing set $\mathcal A_\heeg$, let $\p_i$ denote the union of the midpoints of the arcs in  $T_i$ and set
\[\mathcal A_\heeg=\coprod_{i\in\Z/3}\p_i\times e_i.\]
The triangle formulas below apply to a fixed decorated trisection for which the pointed triple diagram represents the corresponding map. Changing the dividing set can change the link-cobordism map; compare the distinct band decorations in \cite[Section~6.1]{Zemke-1}. Nevertheless, it is shown in \cite{Ef-TQFT} that a substantial part of the information encoded in such cobordism maps remains invariant under decoration change.

For a fixed underlying pair $(X,S)$, the pointed trisection diagram is not unique. Stabilization, destabilization, perturbation and deperturbation are among the moves used to compare such diagrams; see \cite[Sections~3--4]{AE-2} for the related diagram and Cerf moves in tangle Floer theory. {Diagram moves preserving the decorated, multi-based cobordism induce the usual identifications of its map. Moves changing the number of basepoints require the quasi-stabilization maps discussed below; they are not silently identified with isomorphisms of the collapsed complexes.} When comparing decorated cobordisms one must also keep track of the dividing set and the boundary basepoints. We make no identification of maps associated with different dividing sets in this paper, and refer the reader to \cite{Ef-TQFT} for a complete discussion.

\subsection{Heegaard-Floer TQFT}
We use the two-variable specialization of the multi-based link complex in \cite[Remark~2.14 and Section~3.3]{Zemke-1}. Let $\link = (L, \boldsymbol{w}, \boldsymbol{z})$ be a pointed link in a $3$-manifold $Y$, and $\sps \in \SpinC (Y)$ be a $\spinc$ structure on $Y$. Let $\hee = (\Sigma, \boldsymbol{\alpha},  \boldsymbol{\beta}, \boldsymbol{w}, \boldsymbol{z})$ be a corresponding strongly $\sps$-admissible Heegaard diagram and  $\chai^-_{\hee, \sps} = \cfm{}(Y, \link, \sps)$ denote the link Floer chain complex associated with $\hee$ and $\sps$, which is freely generated by the set $\gen (\hee, \sps) \subset \mathbb{T}_\alpha \cap \mathbb{T}_\beta$ of generators representing the $\SpinC$ structure $\sps$, over the polynomial ring $\ftwo [\uu, \vv]$. {Throughout, $\ftwo=\mathbb F_2$. We retain the superscript $-$ for this fully unblocked two-variable complex, and reserve $\circ$ for its $z$-blocked specialization.} The differential $d^-_{\hee, \sps}$ is defined as follows. Given a pair of generators $\gene{x}, \gene{y} \in \gen (\hee, \sps)$ let $\pi^j_2(\gene{x}, \gene{y})$ denote the set of homotopy classes of Whitney disks connecting $\gene{x}$ to $\gene{y}$ with Maslov index $j$. Define 
$$d^-_{\hee, \sps}(\gene{x}) = \sum_{\gene{y} \in \gen (\hee, \sps)} \sum_{\phi \in \pi^1_2(\gene{x}, \gene{y})} \# (\mathcal{M}(\phi)/\mathbb{R}) \cdot \uu^{n_{\boldsymbol{w}}(\phi)} \cdot \vv ^{n_{\boldsymbol{z}}(\phi)} \cdot \gene{y}$$
where $\mathcal{M}(\phi)$ denotes the moduli space of $J$-holomorphic representatives of $\phi$ for a generic path of almost complex structures $J$ which is suppressed in the notation, and 
\[n_{\boldsymbol{w}}(\phi) = \sum_{w \in \boldsymbol{w}} n_w(\phi)\] 
where $n_w(\phi)$ is the coefficient of the region containing $w$ in $\doma(\phi)$. Define  $n_{\boldsymbol{z}}(\phi)$ similarly. We denote the homology of this chain complex by $\hfm{} (Y, \link, \sps)$ which is an $\ftwo [\uu, \vv]$-module. The latter homology group is an invariant for the pointed link $(Y, \link, \sps)$. Let $\hfm{}(Y, \link)$ denote the direct sum of these homology groups for all $\SpinC$ structures:
\[\hfm{}(Y, \link)=\bigoplus_{\sps \in \SpinC (Y)}\hfm{} (Y, \link, \sps).\]
Note that different admissible Heegaard diagrams may be needed for computing the latter group. 

Following the relative grading formulas in \cite[Equation~(2)]{Zemke-1}, we can define \emph{relative homological gradings} using the basepoints in $\boldsymbol{w}$ and $\boldsymbol{z}$ in the following way. The homological bidegrees of $\uu$ and $\vv$ are defined to be $(-2, 0)$ and $(0, -2)$ respectively. Together with the bigradings of the generators, this gives a bigrading on the chain complex $(\chai^-_{\hee, \sps}, d^-_{\hee, \sps})$. The differential drops the bigrading by $(1, 1)$. Note that the  bigrading on the link Floer homology is usually defined in a slightly different (but equivalent) way. 
All grids considered below represent links in $S^3$, where we use the absolute normalizations specified in Section~\ref{subsec:floer-ideals}. 
Define
$$\cfp{} (Y, \link, \sps) := \cfm{}(Y, \link, \sps) \otimes_{\ftwo [\uu, \vv]} \ftwo [\uu, \vv]/\langle \vv \rangle$$
$$\cfp{} (Y, \link) := \bigoplus_{\sps \in \SpinC(Y)} \cfp{} (Y, \link, \sps)$$
and denote the homology of the chain complex $\cfp{} (Y, \link, \sps)$ by $\hfp{} (Y, \link, \sps)$. 

Setting \(v=0\) is a chain-level specialization: precisely the disks
with \(n_{\z}=0\) remain. It is not, in general, interchangeable with taking
homology and then tensoring. The analogous \(w\)-blocked complex is obtained
by setting \(u=0\). A cycle or a nonzero homology class in one of these
specializations need not have the same property in the fully unblocked complex.
All polynomial-valued triangle maps are taken in a fixed Spin$^c$ summand;
an unrestricted sum over Spin$^c$ structures is not asserted to be polynomial.

For the link Floer complex and its naturality, see \cite[Sections~3.3--3.4]{Zemke-1}; compare the tangle Floer naturality result \cite[Theorem~1.1]{AE-2}.

Let $\heeg = (\Sigma, \boldsymbol{\alpha}, \boldsymbol{\beta}, \boldsymbol{\gamma}, \bf{w}, \bf{z})$ be a trisection diagram equipped with compatible dividing data $\mathcal A_\heeg$, and let $\spt \in \SpinC(X_\heeg)$ be such that $\heeg$ is $\spt$-admissible. Let $(Y_i,\link_i)$ denote its boundary pointed links and put $\spt_i=\spt|_{Y_i}$. The holomorphic polygon construction is given in \cite[Section~2.5]{AE-2}, and the $\SpinC$-refined triangle map in \cite[Section~8.3]{OS-3m1}. Under the fixed decorated identification assumed here, the associated triangle maps on the link Floer complexes are
\[
f^\star_{\heeg,\mathcal A_\heeg,\spt}:\chai^\star_{\hee_0,\spt_0}\otimes\chai^\star_{\hee_1,\spt_1}
\longrightarrow\chai^\star_{\hee_2,\spt_2},\qquad \star=-,\circ,
\]
where tensor products and maps are taken over $\ftwo[\uu,\vv]$ for $\star=-$ and over $\ftwo[\uu]$ for $\star=\circ$. For generators $\gene{x}_0$ and $\gene{x}_1$, the formulas, with the output generator also summed over, are
\begin{align*}
f^-_{\heeg,\mathcal A_\heeg,\spt}(\gene{x}_0\otimes\gene{x}_1)
&=\sum_{\gene{x}_2}\sum_{\substack{\Delta\in\pi^0_2(\gene{x}_0,\gene{x}_1,\gene{x}_2)\\ \mathfrak{s}(\Delta)=\spt}}
\#\mathcal M(\Delta)\,\uu^{n_{\w}(\Delta)}\vv^{n_{\z}(\Delta)}\gene{x}_2,\\
f^\circ_{\heeg,\mathcal A_\heeg,\spt}(\gene{x}_0\otimes\gene{x}_1)
&=\sum_{\gene{x}_2}\sum_{\substack{\Delta\in\pi^0_2(\gene{x}_0,\gene{x}_1,\gene{x}_2)\\ \mathfrak{s}(\Delta)=\spt,\ n_{\z}(\Delta)=0}}
\#\mathcal M(\Delta)\,\uu^{n_{\w}(\Delta)}\gene{x}_2.
\end{align*}
Here $\pi^0_2(\gene{x}_0,\gene{x}_1,\gene{x}_2)$ is the set of Maslov-index-zero Whitney triangle classes, $\mathfrak{s}(\Delta)=\mathfrak{s}_{\w}(\Delta)$ is the $w$-based $\SpinC$ structure represented by $\Delta$, and $\mathcal M(\Delta)$ is the moduli space of holomorphic representatives. Composing the natural map from the tensor product of homologies to the homology of the tensor-product complex with the induced triangle map gives the bilinear homology map 
\[
\begin{gathered}
F^\star_{\heeg,\mathcal A_\heeg,\spt}:\hfs{}(Y_0,\link_0,\spt_0)
\otimes\hfs{}(Y_1,\link_1,\spt_1)
\longrightarrow\hfs{}(Y_2,\link_2,\spt_2),\qquad\star=-,\circ.
\end{gathered}
\]
All these maps retain the choice of dividing data in their notation; our later combinatorial formula concerns the $\circ$-triangle map for the fixed decorated trisection.

\subsection{\Goo domains for pointed trisection diagrams} \label{Goo}
Let $\heeg = (\Sigma, \boldsymbol{\lambda_0}, \boldsymbol{\lambda_1}, \boldsymbol{\lambda_2}, \bf{w}, \bf{z})$ be a pointed trisection diagram. Fix generators $\gene{x}_i \in \mathbb{T}_{\lambda_i} \cap \mathbb{T}_{\lambda_{i+1}}$ for $i \in \Z/3$ and a Whitney triangle class $\Delta \in \pi^0_2(\gene{x}_0, \gene{x}_1, \gene{x}_2)$. Then, there is a domain $\doma(\Delta)$ which is a linear combination of 
the connected components in $\Sigma - \boldsymbol{\lambda_0} -\boldsymbol{\lambda_1} - \boldsymbol{\lambda_2}$ (called basic domains) representing the image of $\Delta$ on the surface in the cylindrical reformulation of Heegaard Floer theory \cite{Robert-cylindrical}. By positivity of holomorphic domains, in order for $\Delta$ to correspond to a non-empty moduli space $\mathcal{M}(\Delta)$ of $J$-holomorphic representatives, the coefficients of $\doma(\Delta)$ must be non-negative (provided that the path $J$ of almost complex structures is standard near the marked points). Additionally, the coefficients around every intersection must satisfy a specific equation, defined as follows (see \cite[Section~4]{Sarkar-index}).

Let $(\lambda,\lambda') \in \boldsymbol{\lambda}_i\times\boldsymbol{\lambda}_{i+1}$ with $i \in \Z/3$
and fix $P \in \lambda \cap \lambda'$.
Let $U_P$ be an open neighborhood of $P$ in $\Sigma$ such that $V_P = U_P \cap \lambda$ and $W_P = U_P \cap \lambda'$ are homeomorphic to $\mathbb{R}$ and $U_P \cap \lambda_{i+2} = \varnothing$. We can define a Cartesian coordinate system using an orientation preserving diffeomorphism $\eta: U_P \rightarrow \mathbb{R}^2$ such that $\eta(V_P)$ and $\eta(W_P)$ are $x$-axis and $y$-axis respectively. Then, basic domains around $P$ are in correspondence with four quadrants of $\mathbb{R}^2$. Denote the coefficient of the domain corresponding to the $j$-th quadrant by $c_j(P)$. The point measure below is Sarkar's average of the four local multiplicities \cite[Section~2]{Sarkar-index}. We define
\begin{align}\label{eq:local-quantities}
		&E(P) = \sum_{j=1}^4(-1)^{j-1}c_j(P),
		&&\mu(P) = \frac{1}{4}\sum_{j=1}^4c_j(P) 
		&&\text{and}
		&&\mu_{\gene{x}_i}(\doma) = \sum_{P \in \gene{x}_i} \mu(P).
\end{align} 
The point measure and $E(P)$ are independent of these choices. 

The alternating corner condition is $E(P)=0$ away from the generators and has value $+1$ at every intersection point which belongs to one of the three corner generators.

\begin{defi}
A domain $\doma=\doma(\Delta)$ is called \emph{admissible} if $\doma\geq0$, the corresponding triangle class $\Delta$ has Maslov index zero, and $n_{\boldsymbol{z}}(\Delta)=0$.  Equivalently, it is a positive domain of an index-zero Whitney triangle class that is eligible to contribute to the $z$-blocked triangle map, before its holomorphic count is taken. 
\end{defi}

The above use of ``admissible domain'' is distinct from weak or strong admissibility of a Heegaard diagram.

\section{Grid triples for cobordisms in $\CPbar$}
Let \(\Sigma=\mathbb R^2/\mathbb Z^2\), with its standard complex orientation,
and let \(\alpha,\beta,\gamma\) be horizontal, vertical, and slope \(-1\)
curves. We use the ordered triangle map
\[
C(\alpha,\beta)\otimes C(\beta,\gamma)\longrightarrow C(\alpha,\gamma).
\]
Each pairwise three-manifold is \(S^3\). With the standard cobordism orientation,
the associated four-manifold is
\[
X_{\heeg}=\CPbar\setminus(B_0^4\sqcup B_1^4\sqcup B_2^4),
\qquad \sigma(X_{\heeg})=-1.
\]
Indeed, after a cyclic relabeling this is the genus-one triple subordinate
to a \(-1\)-framed unknot, as in \cite[Section~6]{OS-triangles-revision}.
Here \(\CPbar\) denotes the projective plane with the opposite of its complex
orientation. Reversing the four-manifold orientation gives the positive
projective-plane trisection picture, but also changes the ordered boundary
and Floer conventions; it is not the triangle map displayed above.

We replace each of \(\alpha,\beta,\gamma\) by \(n\) parallel copies.
These redundant compressing disks do not change the underlying handlebodies
or the four-manifold. They provide the additional bridge arcs and basepoints
needed for the grid links. Thus the same oriented identification of
\(X_{\heeg}\) applies to every triple grid below.

\subsection{Triple grid diagrams}

\begin{figure}[t]
        \def\svgwidth{0.3\textwidth}
        \centering
        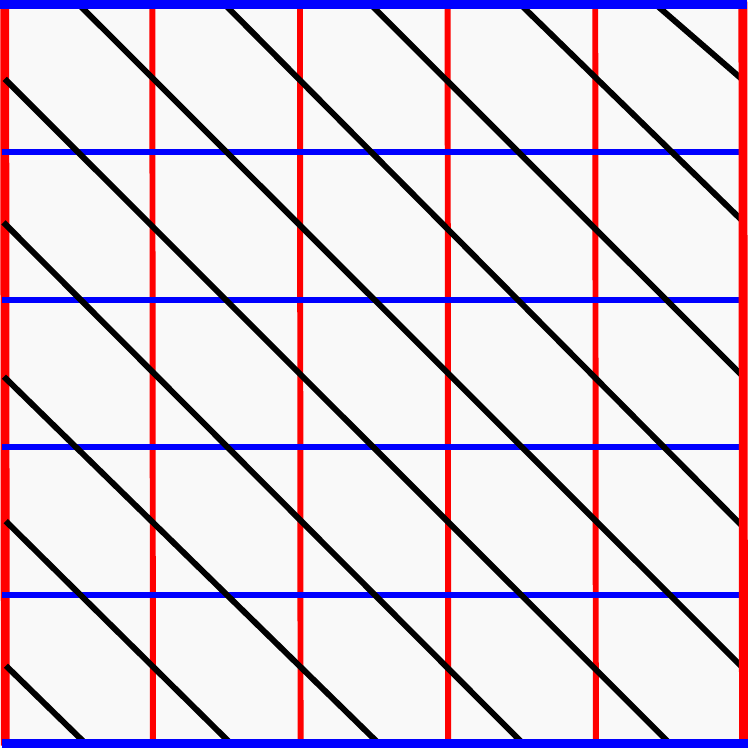
        \caption{The $n\times n$ grid and the $2n$ black diagonal curves representing $\boldsymbol{\gamma}$. The horizontal blue lines form $\boldsymbol{\alpha}$ and the vertical red lines form $\boldsymbol{\beta}$.}\label{fig:triple-grid}
   \end{figure}

Using a Cartesian coordinate system on the universal cover of $\Sigma$, we represent the aforementioned diagrams in a more concrete way.
Let $L_0$, $L_1$ and $L_2$ be three families of lines on the plane defined as follows:
\begin{equation*}
\begin{split}
L_0 &= \{ (x,y) \in \R^2 : y \in \Z \} \\
L_1 &= \{ (x,y) \in \R^2 : x \in \Z \} \\
L_2 &= \{ (x,y) \in \R^2 : x+y=\frac{2k+1}{2},\  k \in \Z \}
\end{split}
\end{equation*}
Consider the torus $\Sigma = \R^2/(n\Z)^2$ and its covering map $\eta : \R^2 \rightarrow \Sigma$. We can define $\boldsymbol{\alpha}$, $\boldsymbol{\beta}$ and $\boldsymbol{\gamma}$ using sets $L_i'=\eta(L_i)$ with $i \in \Z/3$ in the following way. Using the square $(0,n] \times (0,n]$ as a fundamental domain for $\Sigma$, we define coordinates for points of $\Sigma$. In this square $L_0'$ and $L_1'$ form an $n \times n$ grid with $n^2$ unit squares and $L_2'$ forms $2n$ diagonal segments dividing each unit square into three parts, a hexagon and two right angled triangles (see Figure~\ref{fig:triple-grid}). In what follows, whenever we define an object in $\Sigma$ using this coordinate system, we extend the same definition to any object in the preimage of the original object under the map $\eta$.

\begin{defi}
A \emph{triple grid diagram} of \emph{grid number} $n$ is a diagram representing the square $(0,n] \times (0,n]$ accompanied with $L_0'$, $L_1'$ and $L_2'$ as defined above.
\end{defi}

 In this construction all the horizontal, vertical and diagonal lines appearing in $L_0'$, $L_1'$ and $L_2'$ are closed loops on $\Sigma$, and there are $n$ curves of each kind. We call the horizontal, vertical and diagonal curves $\alpha$-curves, $\beta$-curves and $\gamma$-curves, respectively, and color them blue, red and black respectively. Denote the horizontal curve with the equation $y = i$ by $\alpha_i$, and the vertical curve with the equation $x= i$ by $\beta_i$. Note that any $\gamma$-curve on $\Sigma$ is represented by a pair of diagonal segments on the planar diagram. Therefore, we denote the union of diagonal segments with equations 
 \[y = -x + \frac{2i+1}{2}\quad\text{and}\quad y = -x + n + \frac{2i+1}{2}\] 
 by $\gamma_i$. Therefore, given a triple grid diagram, one can construct $(\Sigma, \boldsymbol{\alpha}, \boldsymbol{\beta}, \boldsymbol{\gamma})$ where 
 \[\boldsymbol{\alpha} = \{ \alpha_1, \dots, \alpha_n \},\quad \boldsymbol{\beta} = \{ \beta_1, \dots, \beta_n \}\quad\text{and}\quad \boldsymbol{\gamma} = \{ \gamma_1, \dots, \gamma_n \}\] 
 as above. All indices for curves in $\alphas,\betas$ and $\gammas$ will be understood modulo $n$, i.e. as elements in $\Z/n\Z$.

For $1\leq i,j \leq n$, the $(i,j)$-square of $\Sigma$ is defined as the open square $(i-1,i) \times (j-1, j)$ within $\Sigma$ and $U_{i,j}$, $H_{i,j}$ and $L_{i,j}$ are respectively the upper triangle, the hexagon and the lower triangle inside the $(i,j)$-square of $\Sigma$. Define the set of fundamental domains by
$$\boldsymbol{\Omega} = \{ U_{i,j}, H_{i,j}, L_{i,j} : 1 \leq i, j \leq n \}.$$
The \emph{diagonal value} of the $(i,j)$-square is defined to be the remainder of $i+j-1$ modulo $n$. We also define the $i$th column of $\Sigma$ as the set of $(i, j)$-squares with $1\leq j \leq n$ and similarly, the $j$th row of $\Sigma$ as the set of $(i, j)$-squares with $1\leq i \leq n$.  We also denote the intersection points of the curves in the following way:
$$A_{i,j} = (i, j), \ \ \ B_{i,j} = \Big(i, j -\frac{1}{2}\Big), \ \ \ C_{i,j} = \Big(i- \frac{1}{2}, j \Big)\qquad\text{for}\ \ 1\leq i,j\leq n.$$
Thus $A_{i,j}\in\alpha_{j}\cap\beta_{i}$, $B_{i,j}$ is the intersection of the vertical curve $\beta_{i}$ with $\gamma_{i+j-1}$, and $C_{i,j}$ is the intersection of the horizontal curve $\alpha_{j}$ with $\gamma_{i+j-1}$.  All indices are understood modulo $n$.
	
Let $S_n$ be the symmetric group of the set $\{1,2, \dots, n\}$.
Given permutations  $\sigma_{\w}, \sigma_{\z} \in S_n$, we define two sets of marked points on $\Sigma$:
$$\w = \{ w_1, w_2, \dots, w_{n} \},\ \z= \{ z_1, z_2, \dots, z_{n} \}$$
 We require $\sigma_{\w}(i)\ne\sigma_{\z}(i)$ for every $i$, so that the two sets of basepoints are disjoint. We place each $w_i$ at the center of $(i, \sigma_{\w}(i))$-square and each $z_i$ at the center of $(i, \sigma_{\z}(i))$-square, for $1\leq i \leq n$.   Moreover, we define a diagonal value function 
$$d: {\w \cup \z} \rightarrow \{ 1, 2, \dots, n\}=\Z/n\Z$$
that assigns to each marked point the diagonal value of the square containing it. Therefore, we obtain
$$d(w_i) = i+\sigma_{\w}(i)-1\ (\mathrm{mod} \ n), \ \ d(z_i) = i + \sigma_{\z}(i)-1 \ (\mathrm{mod}\ n).$$ Denote the grid diagrams obtained by removing the diagonal curves, the horizontal curves and the vertical curves by $\heeg_0$, $\heeg_1$ and $\heeg_2$ respectively.  Therefore, these grid diagrams are formed by the pairs $(\boldsymbol{\alpha},\boldsymbol{\beta})$,
$(\boldsymbol{\beta},\boldsymbol{\gamma})$ and
$(\boldsymbol{\alpha},\boldsymbol{\gamma})$, respectively.

\begin{defi}
Let $\sigma_{\w}$ and $\sigma_{\z}$ be a pair of permutations in $S_n$ and define marked points $\w$ and $\z$ as above. We call the marked points $\w$ and $\z$ \emph{appropriate} if both of $\heeg_1$ and $\heeg_2$ have exactly one marked point of each kind in every row and every column.
\end{defi}

Equivalently, we need $(\Sigma, \boldsymbol{\alpha}, \boldsymbol{\beta}, \boldsymbol{\gamma}, {\w }, {\z})$ to be a pointed trisection diagram, or that each connected component of $\Sigma - \boldsymbol{\alpha}$, $\Sigma - \boldsymbol{\beta}$ and $\Sigma - \boldsymbol{\gamma}$ contains precisely one marking from $\w$ and one marking from $\z$.

Appropriateness is equivalent to requiring that both maps
\(i\mapsto i+\sigma_{\w}(i)\) and \(i\mapsto i+\sigma_{\z}(i)\)
are permutations modulo \(n\). In particular \(n\) is odd: summing either
permutation condition gives \(n\mid n(n-1)/2\). Since disjoint central
markings exclude \(n=1\), the centrally marked grids below have \(n\ge3\).

\begin{figure}[t]
	\def\svgwidth{0.9\textwidth}
	\centering
	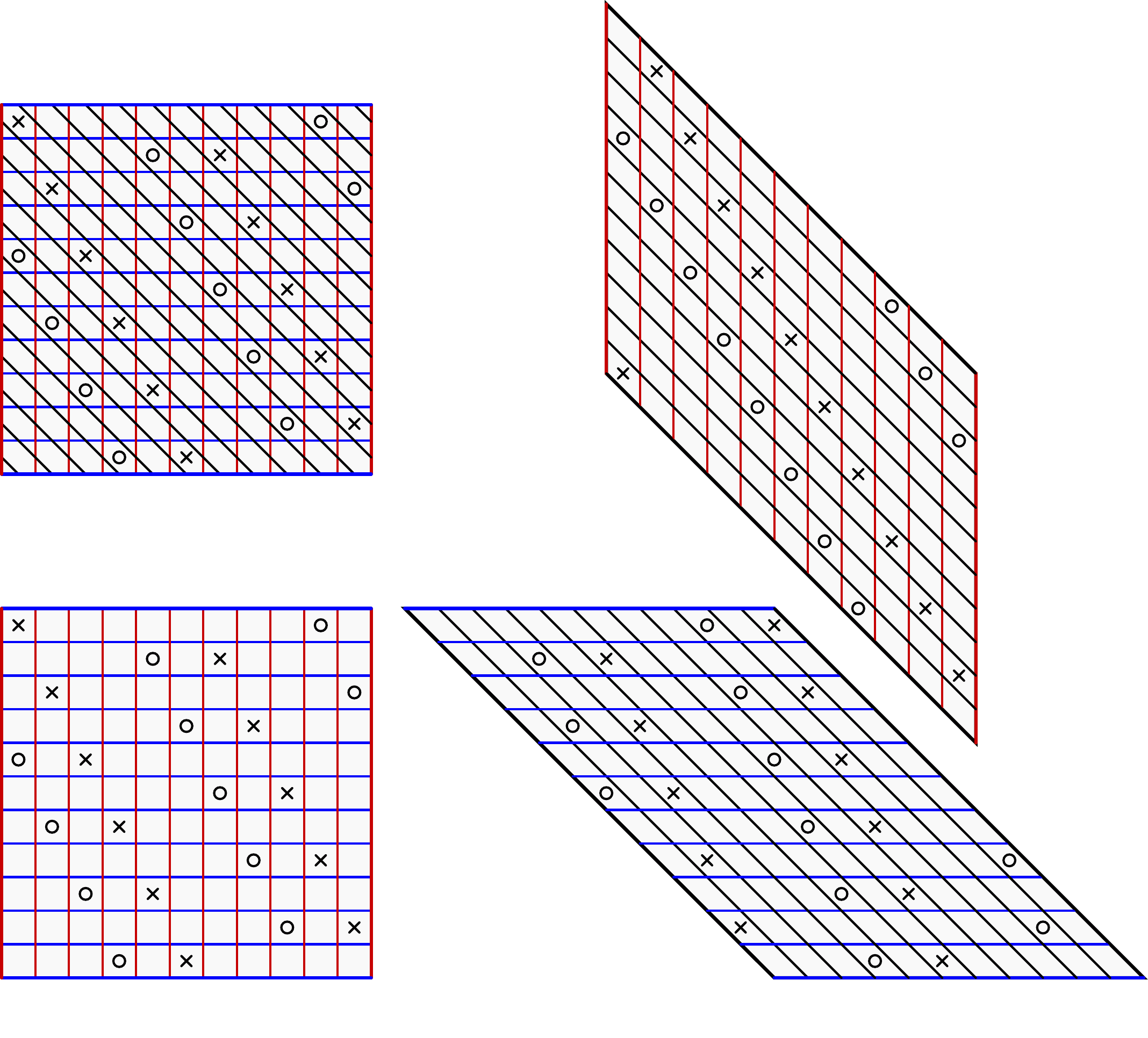
	\caption{Three grid diagrams $\heeg_0$, $\heeg_1$ and $\heeg_2$ are constructed from the triple grid diagram $\heeg$ by removing the black, blue and red line segments, respectively, and changing the fundamental domain, when needed (parts (c) and (d)).}\label{fig:G0G1G2}
\end{figure}

We can construct the diagrams $\heeg_1$ and $\heeg_2$ by changing the fundamental domain of $\Sigma$. Let  $\mathbb{D}_1$ be the parallelogram whose sides are located on lines with equations 
\begin{align*}
	&x = 0, &&x = n, &&y = -x + n + \frac{1}{2} &&\text{and} 
	&&y = -x + \frac{1}{2}.
\end{align*}	 
Similarly,  define
$\mathbb{D}_2$ as the parallelogram whose sides are segments located on lines with equations 
\begin{align*}
&y = 0, &&y = n, &&y = -x + n + \frac{1}{2} &\text{and} 
&&y = -x + 2n + \frac{1}{2}.
\end{align*}
 Then $\mathbb{D}_1$ and $\mathbb{D}_2$ are fundamental domains for $\Sigma$ and diagonal curves are not divided in these diagrams as shown in parts (c) and (d) of Figure~\ref{fig:G0G1G2}. 

 Every triple grid diagram with appropriate marked points gives a pointed trisection diagram for a link cobordism in $\CPbar$:
$$\heeg = (\Sigma, \boldsymbol{\alpha}, \boldsymbol{\beta}, \boldsymbol{\gamma}, {\w } = \{ w_1, w_2, \dots , w_{n} \}, {\z} = \{ z_1, z_2, \dots , z_{n} \})$$

We abuse the notation and identify $\heeg$ with the corresponding triple grid diagram.

\subsection{The genus and the homology class of the decorated surface}
We next record the homology class and the Euler characteristic of the decorated surface \(S_\heeg\) determined by a triple grid \(\heeg\).

\begin{figure}[t]
\def\svgwidth{0.3\textwidth}
\centering
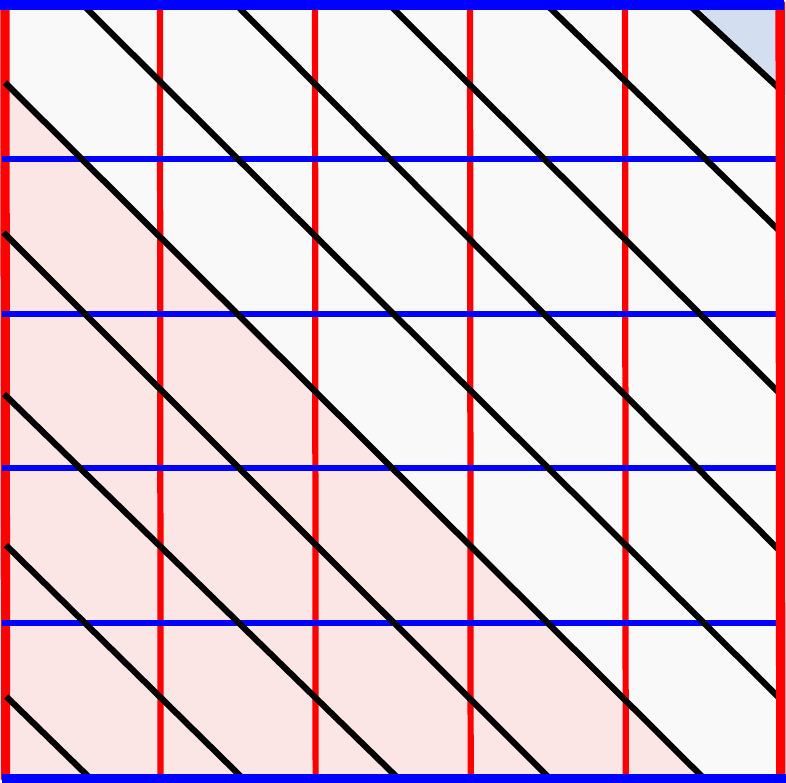
\caption{The reference triply periodic domain \(\Pcal_0\). The red and
	blue regions have coefficients \(1\) and \(-1\), respectively.}
\label{fig:difference-domain}
\end{figure}

Define a triply periodic domain \(\Pcal_0\) by setting
\begin{align*}
	&\Pcal_0=\sum_{1\le i,j\le n}(l_{i,j}L_{i,j}+h_{i,j}H_{i,j}+u_{i,j}U_{i,j}),\qquad\text{where}\\
	&l_{i,j}(\Pcal_0)=\mathbf1_{\{ i+j\le n+1\}},\qquad
	h_{i,j}(\Pcal_0)=\mathbf1_{\{i+j\le n\}}\qquad\text{and}\qquad
	u_{i,j}(\Pcal_0)=\mathbf1_{\{i+j< n\}}-\mathbf1_{\{(i,j)=(n,n)\}}.
\end{align*}
This triply periodic domain,  illustrated in Figure~\ref{fig:difference-domain}, is in fact the primitive periodic domain in the genus-one \(-1\)-framed-unknot triple once the redundant $\alpha$-,$\beta$- and $\gamma$-curves are removed and just $\alpha_n,\beta_n$ and $\gamma_{n-1}$ are preserved. We fix \(H=[\Pcal_0]\in H_2(X_{\heeg};\mathbb Z)\) as the homology class obtained by attaching disks to the boundary of $\Pcal_0$ (i.e. to $\alpha_n,\beta_n$ and $\gamma_{n-1}$) in $X_\heeg$, so that \(H^2=-1\).
 Let \(H_{\rm rel}\) denote  the image of $H$ in relative homology.

\begin{prop}\label{prop:shadow-degree}
If \([S_{\heeg}]=d(\heeg)H_{\rm rel}\), then
\[
d(\heeg)=\sum_{w\in\w}n_{\w}(\Pcal_0)-\sum_{z\in\z}n_{\z}(\Pcal_0).
\]
\end{prop}
\begin{proof}
In the standard orientation the surface sheets over the parameter triangle
are positive at \(z\) and negative at \(w\). The caps for \(\Pcal_0\) may be
chosen disjoint from the bridge arcs. Therefore
\[
H\cdot[S_{\heeg}]=n_{\z}(\Pcal_0)-n_{\w}(\Pcal_0).
\]
The left side is \(-d(\heeg)\), since \(H\cdot H_{\rm rel}=-1\).
This proves the formula, and agrees with
\cite[Equation~(3.2)]{Zemke-gradings-revision}.
\end{proof}

\begin{prop}\label{prop:nullhomologous-grid-surface}
For the appropriate central markings used in this paper, the surface is
relatively nullhomologous:
\[
[S_{\heeg}]=0\quad\text{in}\quad
H_2(X_{\heeg},\partial X_{\heeg};\mathbb Z).
\]
\end{prop}
\begin{proof}
The periodic domain \(\Pcal_0\)  has hexagonal coefficients
\[
h_{i,j}(\Pcal_0)=\begin{cases}1,&i+j\le n,\\0,&i+j>n.\end{cases},\qquad\forall\ \ i,j\in\{1,\ldots,n\}.
\]
Every marking is in a hexagon. For either marking permutation \(\sigma\),
let \(r_i\in\{1,\ldots,n\}\) be congruent to \(i+\sigma(i)\).
Appropriateness makes the \(r_i\) a permutation, and hence
\begin{align*}
&\#\{i:i+\sigma(i)>n\}
=\sum_{i=1}^n\frac{i+\sigma(i)-r_i}{n}
=\frac{n+1}{2}\\
\Rightarrow\qquad&
\sum_{i=1}^n h_{i,\sigma(i)}(\Pcal_0)=n-\#\{i:i+\sigma(i)>n\}=\frac{n-1}2.
\end{align*}
It follows from the above computation, which is the same for $\sigma_{\w}$ and $\sigma_{\z}$, and the preceding proposition, that 
\(d(\heeg)=0\).
\end{proof}

\begin{prop}\label{prop:general-genus}
Let \(S_\heeg\) be the oriented surface associated to a pointed grid triple $\heeg$ of grid number \(n\). Let \(b\) be the total number of boundary components,
\(c\) the number of connected components of \(S_\heeg\), and \(g_{\rm tot}\)
the sum of their genera. Then
\[
\chi(S_\heeg)=-n,\qquad g_{\rm tot}=c+\frac{n-b}{2}.
\]
In particular, if \(S_\heeg\) is connected then
\[
g(S_\heeg)=\frac{n+2-b}{2},
\]
and if all three boundary links are knots then \(g(S_\heeg)=(n-1)/2\).
\end{prop}
\begin{proof}
The surface has one disk \(\{v\}\times\triangle\) for each of the \(2n\)
basepoints and one band for each bridge arc in each of the three curve
families, hence \(3n\) bands. Therefore \(\chi(S_\heeg)=2n-3n=-n\).
The formula \(\chi=2c-2g_{\rm tot}-b\) gives the genus statements.
If any boundary link is a knot, its boundary circle meets every disk in the
corresponding two-colored graph, so the whole surface is connected.
\end{proof}

\subsection{\Nic domains for a triple grid diagram}
Let \(\heeg\) be a triple grid of grid number \(n\) with appropriate markings.

\begin{defi}
Let \(\eta:\mathbb R^2\to\Sigma\) be the covering map. A \emph{nice triangle}
$\tria$ is the push-forward, with multiplicities, of a closed Euclidean triangle
$\widetilde \tria$ whose sides lie on horizontal, vertical, and diagonal lifts of attaching curves. The equal shorter sides of the right isosceles triangle $\widetilde \tria$ have length \(k-\tfrac12\), where \(|\tria|=k\ge1\) is a positive
integer. The triangle is \emph{empty} if its domain has \(n_{\z}(\tria)=0\). Its {\emph{corners}} are the projected vertices.
\end{defi}

\begin{defi}
Let \(\psi\in\pi_2(\x_0,\x_1,\x_2)\). Its domain is \emph{nice} if
\[
\mathcal D(\psi)=\mathcal T_1+\cdots+\mathcal T_n,
\]
where  \(\mathcal T_i\) are nice triangles with pairwise distinct
\(\alpha\)-, \(\beta\)-, and \(\gamma\)-boundary curves, and their corners
are the specified generators, while $\mu(\psi)=0$.
It is an \emph{empty nice} domain if every summand is empty.
\end{defi}

\begin{thm}\label{mai}
For the split complex structure (specified in Section~\ref{holo}) every positive
index-zero triangle class $\psi$ with a holomorphic representative has 
a nice domain. If \(n_{\z}(\psi)=0\), each triangular sheet is embedded and has side length at most \(n-\tfrac12\); in particular it has an empty nice domain.
Conversely every nice class has mod-\(2\) count one in this model.
\end{thm}

The theorem is proved in Section~\ref{holo}.
For a nice domain, following the boundary of a sheet encounters one
\(\alpha\beta\)-corner, one \(\beta\gamma\)-corner, and one
\(\alpha\gamma\)-corner before returning to the starting point.
This sheetwise rule remains meaningful when projected domains overlap;
such overlap is not the same as intersection with the symmetric-product diagonal.

\subsection{$\SpinC$ structures for triple grid diagrams}
\label{sec:spinc-triple-grid}
Write \(X=X_{\heeg}\). For generators \(\x_0,\x_1,\x_2\) in the three
pairwise diagrams, the Whitney-triangle construction associates the $\SpinC$ structures
\(\spinc_{\w}(\psi)\) and \(\spinc_{\z}(\psi)\) in $\SpinC(X)$; see
\cite{OS-3m1,AE-2}. We label cobordism maps by \(\spinc_{\w}\).
The identity
\[
\spinc_{\w}(\psi)-\spinc_{\z}(\psi)=\PD[S_{\heeg}]
\]
is \cite[Lemma~3.9]{Zemke-gradings-revision}. By
Proposition~\ref{prop:nullhomologous-grid-surface} these two structures agree
for the grids in this paper; we write their common value as \(\spinc(\psi)\).
The boundary structures on \(S^3\) are unique, and \(\SpinC(X)\) is an
affine space over \(H^2(X;\mathbb Z)\cong\mathbb Z\).

For later reference, let \(\psi_0\) be the product of the small upper triangles
\(U_{i,j}\) in the squares containing the \(\z\)-markings. Appropriateness
ensures that its corners are three grid states. The triangles are translates
of one another, so the product of their conformal parametrizations avoids
the fat diagonal and has Maslov index zero. We denote the $\SpinC$ class represented by $\psi_0$ by $\spinct_0$.

\begin{defi}
For a formal domain
\[
\mathcal D=\sum_{i,j}\bigl(u_{i,j}U_{i,j}+h_{i,j}H_{i,j}+l_{i,j}L_{i,j}\bigr),
\]
define its \emph{difference number} by
\(\delta(\mathcal D)=\sum_{i,j}(l_{i,j}-u_{i,j})\).
For a triangle class put \(\delta(\psi)=\delta(\mathcal D(\psi))\).
\end{defi}

\begin{prop}\label{prop:difference-number}
If a triply periodic domain \(\Pcal\) represents \(mH\), then
\[
\delta(\Pcal)=2nm.
\]
\end{prop}
\begin{proof}
The displayed coefficients give
\[
\delta(\Pcal_0)=\frac{n(n+1)}2-\frac{(n-2)(n-1)}2+1=2n.
\]
The functional \(\delta\) vanishes on every domain supported on a pairwise
grid: each elementary pairwise parallelogram contains equal total numbers
of upper and lower triangles. This includes pairwise Whitney-disk domains
and doubly periodic domains. It also vanishes on \([\Sigma]\).
Modulo the pairwise periodic relations and the whole-surface domain, the
triply periodic domains identify with \(H_2(X;\mathbb Z)\), generated by
\(\Pcal_0\). Thus \(\Pcal-m\Pcal_0\) is a sum of domains
annihilated by \(\delta\), proving the formula.
If basepoint normalization is required, subtract row or column domains to
make the individual basepoint multiplicities zero. Such normalization does
not change \(\delta\) or the represented class.
\end{proof}

Let \(\Tri(\heeg)\) denote the set of all Whitney-triangle classes,
with their corners allowed to vary, and put \(h=-\PD[H]\in H^2(X;\mathbb Z)\); 
in particular \(\langle h,H\rangle=1\).

\begin{thm}\label{thm:delta-on-triangles}
The function
\[
\widetilde\delta:\Tri(\heeg)\longrightarrow\tfrac12+\mathbb Z,
\qquad \widetilde\delta(\psi)=\frac{-\delta(\psi)}{2n},
\]
depends only on \(\spinc(\psi)\), satisfies
\(\widetilde\delta(\psi_0)=\tfrac12\), and obeys
\[
\spinc(\phi)-\spinc(\psi)
=\bigl(\widetilde\delta(\phi)-\widetilde\delta(\psi)\bigr)h.
\]
Consequently, if \(c_\psi=\langle c_1(\spinc(\psi)),H\rangle\), then
\[
c_\phi-c_\psi=2(\tdelta(\phi)-\tdelta(\psi)).
\]
\end{thm}
\begin{proof}
Triangle classes representing the same Spin$^c$ structure differ by splicing
pairwise Whitney disks; their domains have equal \(\delta\), by the preceding
proof. After corners have been connected by such disks, the remaining
difference is a triply periodic domain \(\Pcal\), and 
 \[\spinc(\psi+\Pcal)-\spinc(\psi)=\PD[\Pcal].\]
If \([\Pcal]=mH\), then \(\PD[\Pcal]=-mh\), while Proposition~\ref{prop:difference-number} gives \(\widetilde\delta(\psi+\Pcal)-\widetilde\delta(\psi)=-m\). This proves the asserted difference formula.
The reference triangle class consists of \(n\) upper triangles, so
\(\delta(\psi_0)=-n\). This proves half-integrality.
Finally 
\[c_1(\spinct+mh)=c_1(\spinct)+2mh,\] 
and evaluation on
\(H\) changes the Chern number by \(2m\).
\end{proof}
An absolute calibration of this affine coordinate is given in
Corollary~\ref{cor:absolute-chern-coordinate} below.

\section{Computation of the link cobordism map}\label{comp}

In this section, we prove that the holomorphic triangle map associated with a grid triple may be computed combinatorially. Let us start with the index computations.

\subsection{Euler measure}\label{subsec:euler-revised}
Let
\[
\doma=
\sum_{i,j\in\mathbb{Z}/n\mathbb{Z}}
\bigl(u_{i,j}U_{i,j}+h_{i,j}H_{i,j}+l_{i,j}L_{i,j}\bigr)
\]
be the domain of a Whitney triangle 
\[
\psi\in\pi_2(\gene{x}_0,\gene{x}_1,\gene{x}_2).
\]
In particular,  its corner function is $1$ at the coordinates of the relevant generator and $0$ at every other pairwise intersection.
For later use, we also set
\[
\begin{aligned}
a_{i,j}
&=u_{i,j}+l_{i+1,j+1}-h_{i,j+1}-h_{i+1,j},\\
b_{i,j}
&=u_{i,j}+l_{i+1,j}-h_{i,j}-h_{i+1,j},\\
c_{i,j}
&=u_{i,j}+l_{i,j+1}-h_{i,j}-h_{i,j+1}.
\end{aligned}
\]
Then $a_{i,j}=1$ exactly when $A_{i,j}\in\gene{x}_0$,
$b_{i,j}=1$ exactly when $B_{i,j}\in\gene{x}_1$, and
$c_{i,j}=1$ exactly when $C_{i,j}\in\gene{x}_2$; all other values are
$0$.

We use the Euler measure convention of \cite{Sarkar-index}.  The upper and
lower triangles have Euler measure $1/4$, while a hexagonal region has
Euler measure $-1/2$.  Therefore
\[
e(\doma)
=
\frac14\sum_{i,j}
\bigl(u_{i,j}+l_{i,j}-2h_{i,j}\bigr).
\]
On the other hand,
\[
E(A_{i,j})
=u_{i,j}+l_{i+1,j+1}-h_{i,j+1}-h_{i+1,j},
\]
and hence, after reindexing the sums,
\[
\sum_{i,j}E(A_{i,j})
=
\sum_{i,j}
\bigl(u_{i,j}+l_{i,j}-2h_{i,j}\bigr).
\]
By the defining corner condition, $E(A_{i,j})=1$ at the $n$
coordinates of $\gene{x}_0$ and is $0$ elsewhere.  We have proved:

\begin{lem}
\label{lem:euler-measure}
Every domain $\doma=\Dcal(\psi)$ associated with a Whitney triangle $\psi$ satisfies
\[
e(\doma)=\frac n4.
\]
\end{lem}

\subsection{Holomorphic representatives}\label{holo}
Equip \(\Sigma=\mathbb C/(n\mathbb Z+in\mathbb Z)\) with its flat complex
structure, and its symmetric product with the induced complex structure.
Away from the fat diagonal, the ordered configuration cover has the product
complex structure. We compute index-zero triangle counts in this split
model. The representatives below are regular, so these counts may be used
in the Floer construction and transported by continuation through the
standard admissible almost-complex structures preserving the diagonal.
Indeed the index formula and positivity prohibit negative-index
holomorphic triangles; disk or sphere bubbling therefore do not give
an extra end in the zero-dimensional triangle count.

\begin{thm}\label{thm:holomorphic-domains-are-nice}
In this split model, every positive index-zero triangle class admitting
a holomorphic representative has a nice domain. If it is $\z$-blocked,
it has an empty nice domain and for every nice triangle $\tria$ appearing in the domain we have $|\tria|\le n$.
\end{thm}
\begin{proof}
Sarkar's index formula \cite[Theorem~4.1]{Sarkar-index} is
\[
\mu(\psi)=\iota(\psi)+2e(\mathcal D)-\frac n2,
\]
where \(\iota\) is intersection with the fat diagonal.
Lemma~\ref{lem:euler-measure} therefore gives \(\iota(\psi)=0\).

The diagonal is a complex divisor and the boundary tori are disjoint from it.
Positivity of intersections implies that a holomorphic representative avoids
the diagonal. The source triangle is simply connected, so the map lifts from $\mathrm{Sym}^n(\Sigma)$ to $\Sigma^n$ and thus splits into \(n\) holomorphic maps to \(\Sigma\). Each has one boundary arc on a curve of each family.

Lift a component to \(\mathbb C\). Its boundary arcs lie on three fixed
straight lines. For a point off these lines, the winding number of its
boundary equals that of the polygon obtained by replacing each boundary
arc by the straight segment with the same endpoints: excursions along a
line do not change that winding number. The argument principle implies
that the component maps with degree one onto the bounded triangle and
has no image outside it. It is therefore a conformal parametrization of
that Euclidean triangle. The generator condition supplies distinct
attaching curves for the \(n\) sheets, The product avoids the diagonal, so Proposition~\ref{prop:sheet-separation} gives the separation condition, proving niceness.

Suppose now \(n_{\z}(\mathcal D)=0\). A right isosceles triangle with the
prescribed directions and leg length \(L<n\) projects injectively, since
both coordinate ranges have length less than \(n\). If instead
\(L\ge n\), its half-integral length is at least \(n+\tfrac12\).
In the row of square centers one half-unit from its horizontal leg,
its interior contains representatives of all \(n\) columns. The unique
\(\z\)-marking in that row is therefore in the triangle, a contradiction.
Thus \(L\le n-\tfrac12\), and every sheet is empty and embedded.
The same argument, with \(w\) in place of \(z\), applies to a \(w\)-blocked
representative.
\end{proof}

\begin{lem}\label{lem:nice-count}
Every nice class has  mod-$2$ holomorphic count one in the
split model. In particular, this holds for every empty nice class.
\end{lem}
\begin{proof}
For each triangular summand choose its lift to \(\mathbb C\).
The Riemann mapping theorem gives a unique conformal map from the common
source triangle to that lift with the three ordered vertices prescribed.
Their product, followed by the symmetric-product quotient, is holomorphic
and has the prescribed domain.

The linearized problem consequently splits into the problems for the
individual planar triangles. Each has index zero and trivial kernel:
an infinitesimal holomorphic deformation tangent to all three straight
sides, with fixed vertices, pulls back to a holomorphic vector field on
the disk tangent to its boundary and vanishing at three boundary points,
and hence is zero. Equivalently this is the standard three-convex-corner
Riemann--Hilbert problem. Index zero now gives surjectivity. Thus the product
representative is regular and isolated.

Each attaching curve occurs on only one sheet. Its oriented boundary
one-chain therefore fixes the corresponding lifted side, including any
winding; hence the domain determines its triangular sheet matching.
Any other representative avoids the diagonal and splits into those same
triangles by the preceding theorem. Uniqueness of their conformal maps
proves that the moduli space has one point.
\end{proof}

Given a $\SpinC$ class of triangles and having fixed the corners $\x_0,\x_1,\x_2$, there are
only finitely many candidate sheets which have index $0$. 
The separation test \eqref{eq:sheet-separation} provides 
us with a more friendly test for computing the Maslov index, and checking whether it is zero for a domain which is a sum of triangles.
Therefore, it gives a more concrete combinatorial algorithm for finding the contributing triangles.

\subsection{Separation of triangular sheets and finiteness}\label{subsec:sheet-separation}
We now give a geometric test for the index-zero condition which uses only the
vertices and side lengths of the triangular sheets. Put
\[
T=\{(x,y)\in\mathbb R^2:x\ge0,\ y\ge0,\ x+y\le1\}.
\]
For a compatible collection of \(n\) nice triangles, choose lifts
\(A_i=(a_i,b_i)\in\mathbb Z^2\) of their right-angle vertices and signed
leg lengths \(\lambda_i\in\tfrac12+\mathbb Z\), so that their lifted
triangles are \(A_i+\lambda_iT\). Positive lengths correspond to triangles
extending upwards and rightwards, and negative lengths to those extending
downwards and leftwards. The affine maps
\(f_i(q)=\eta(A_i+\lambda_iq)\) parametrize the sheets with the same
ordered source corners. Compatibility means that the sheets use distinct
attaching curves in each family and have the prescribed corner generators.
Their product defines a holomorphic Whitney triangle \(\psi\) in the split
model, without any initial assumption on its index.

\begin{prop}\label{prop:sheet-separation}
For each \(i<j\), define the difference triangle
\[
R_{ij}=\operatorname{conv}\{A_i-A_j,\
A_i-A_j+(\lambda_i-\lambda_j,0),\
A_i-A_j+(0,\lambda_i-\lambda_j)\}.
\]
Then \(\mu(\psi)=\iota(\psi)=0\) if and only if
\begin{equation}\label{eq:sheet-separation}
R_{ij}\cap(n\mathbb Z)^2=\varnothing\qquad\text{for every }i<j.
\end{equation}
Equivalently, the map
\[
\coprod_{i=1}^n T\longrightarrow T\times\Sigma,
\qquad(q,i)\longmapsto(q,f_i(q))
\]
is injective. More generally,
\begin{equation}\label{eq:sheet-index}
\mu(\psi)=\iota(\psi)
=2\sum_{\substack{i<j\\\lambda_i\ne\lambda_j}}
\#\bigl(\operatorname{Int}R_{ij}\cap(n\mathbb Z)^2\bigr).
\end{equation}
\end{prop}
\begin{proof}
Two sheets coincide at \(q\) precisely when
\(A_i-A_j+(\lambda_i-\lambda_j)q\in(n\mathbb Z)^2\).
Compatibility excludes coincidences on the boundary: on each source edge
the sheets lie on distinct curves of the corresponding family. If
\(\lambda_i=\lambda_j\), they are distinct translates and never coincide.
Otherwise each interior lattice point of \(R_{ij}\) gives one simple zero
of their difference. Locally the discriminant of the symmetric-product
quotient is the product of the squared pairwise differences. Each such
zero therefore contributes two to \(\iota\), with contributions adding
also at simultaneous collisions of several sheets. This proves the
formula for \(\iota\). Lemma~\ref{lem:euler-measure} and Sarkar's index
formula give \(\mu=\iota\), proving all the assertions.
\end{proof}
The test can alternatively be expressed without enumerating lattice points. For
\(d=\lambda_i-\lambda_j\ne0\), set
\begin{align*}
\varepsilon=\operatorname{sign}(d),\qquad
r_x=[\varepsilon(a_j-a_i)]_n,\qquad
r_y=[\varepsilon(b_j-b_i)]_n,
\end{align*}
using representatives in \(\{1,\ldots,n-1\}\); these are nonzero by
distinctness of the vertical and horizontal curves. A collision occurs
exactly when \(r_x+r_y<|d|\). Equality is excluded by distinctness of the
diagonal curves. Thus the pair satisfies \eqref{eq:sheet-separation}
exactly when \(|d|<r_x+r_y\). This is a finite test for niceness, independent
of a prior Maslov-index computation. Together with Lemma~\ref{lem:nice-count},
it determines precisely which compatible triangular collections contribute.

\begin{cor}\label{cor:finite-triangle-summands}
For fixed corners and a fixed \(\SpinC\) structure, only finitely many
classes contribute to the unblocked triangle map of the grid triple.
\end{cor}
\begin{proof}
The separation test bounds every difference
\(|\lambda_i-\lambda_j|<2n\). Counting elementary upper and lower triangles
in a sheet gives \(\delta(\mathcal T_i)=2\lambda_i\), so
\(\delta(\psi)=2\sum_i\lambda_i\).
Theorem~\ref{thm:delta-on-triangles} fixes this sum in a fixed
\(\SpinC\) structure. Each signed length is therefore bounded. There are
only finitely many half-integral choices and finitely many right-angle
vertices modulo \((n\mathbb Z)^2\), proving the assertion.
\end{proof}

Individual sheets need not project injectively to \(\Sigma\). For example,
enlarge all the small upper triangles of \(\psi_0\) to the common leg
length \(kn+\tfrac12\), with \(k\ge1\). Their corners are unchanged, and
the sheets are still distinct translates, so they satisfy
\eqref{eq:sheet-separation} and have index zero. Their individual projections
are noninjective, but their difference number is \(-n(2kn+1)\), so their
\(\SpinC\) structures vary with \(k\). In contrast, enlarging an upper
and a lower sheet to this common length gives a signed-length difference
\(2kn+1>2n-2\), forcing a collision and positive index. Thus an increase
in domain multiplicities with fixed corners and fixed \(\SpinC\) structure
does not preserve index zero: for such classes the periodic-domain index
formula gives
\[
\mu(\psi')-\mu(\psi)
=2\bigl(n_{\w}(\psi')-n_{\w}(\psi)\bigr)
=2\bigl(n_{\z}(\psi')-n_{\z}(\psi)\bigr);
\]
see \cite[Lemma~6.10]{Zemke-gradings-revision}.
The fixed-summand finiteness above justifies polynomial-valued maps over
\(\mathbb F[\uu,\vv]\); it does not justify an unrestricted sum over all
\(\SpinC\) structures without a suitable completion.

\subsection{Torsion-free Floer ideals and grid diagrams}\label{subsec:floer-ideals}
Set \(\Abb=\mathbb F[u,v]\), with \(\deg u=(-2,0)\) and
\(\deg v=(0,-2)\).
For a knot \(K\subset Y\) and a $\SpinC$ structure $\spinc\in\SpinC(Y)$, let \(C^{(n)}(Y,K,\spinc)\) denote the two-variable knot Floer chain complex associated with $K$ decorated with $n$ pairs of basepoints, constructed from a compatible Heegaard diagram and freely generated over $\Abb$ by the corresponding grid or Heegaard generators.  Let \(C(Y,K,\spinc)=C^{(1)}(Y,K,\spinc)\) denote the \emph{doubly based} chain complex and define
\[
H(Y,K,\spinc)=H_*(C(Y,K,\spinc)),\qquad
\Ical(Y,K,\spinc)=H(Y,K,\spinc)/\operatorname{Tors}_{\Abb}H(Y,K,\spinc).
\]
By \cite[Proposition~3.4]{AE-unknotting}, when $Y=S^3$ this is a rank-one torsion-free module represented by a monomial ideal after specifying its grading shift.

Given a multi-based Heegaard diagram representing a knot $K\subset Y$, a sequence of diagram identifications and quasi-stabilizations gives a chain homotopy equivalence
\begin{align*}
C^{(n)}(Y,K,\spinc)\simeq C(Y,K,\spinc)\otimes_{\Abb}V^{\otimes(n-1)},
\quad V=\langle\theta^+,\theta^-\rangle_{\Abb},
\quad \partial_V=0.
\end{align*}
In the cobordism grading convention,
\(\deg\theta^+=(\tfrac12,-\tfrac12)\) and
\(\deg\theta^-=(-\tfrac12,\tfrac12)\).
In particular, for \(Y=S^3\), the torsion-free rank before reduction is \(2^{n-1}\); in general it is \(2^{n-1}\operatorname{rank}_{\Abb}\Ical(Y,K,\spinc)\).

Let us now assume that $\heeg=(\Sigma,\alphas,\betas,\gammas,\ws,\zs)$ is a Heegaard triple with $|\z|=|\w|=n$, which is $\spinct$-admissible for some $\spinct\in\SpinC(X_\heeg)$. Let $K_0,K_1$ and $K_2$ denote the corresponding boundary links in the three-manifolds $Y_0,Y_1$ and $Y_2$, respectively, and  $\spinc_i=\spinct|_{Y_i}$. Further assume that each $K_i$ is a knot. Then the associated cobordism map 
\begin{align*}
	f_{\heeg,\spinct}:C^{(n)}(Y_0,K_0,\spinc_0)\otimes_{\Abb}
	C^{(n)}(Y_1,K_1,\spinc_1)\to C^{(n)}(Y_2,K_2,\spinc_2)
\end{align*}	   
induces a corresponding cobordism map 
\begin{align}\label{eq:induced-map}
	F_{\heeg,\spinct}:H(Y_0,K_0,\spinc_0)\otimes_{\Abb}
	H(Y_1,K_1,\spinc_1)\to H(Y_2,K_2,\spinc_2),
\end{align}
which takes torsion classes to torsion classes, 
see \cite[Theorem~3.7]{Ef-TQFT}. Here the reduced homology groups are
identified with the common kernels of the $\z$-basepoint actions in the
multi-based homologies. Equivalently, in the tensor identification above
we use the factor $(\theta^+)^{\otimes(n-1)}$ on each boundary. The analogous
reduction uses the $\w$-actions and $(\theta^-)^{\otimes(n-1)}$; we retain
the $\theta^+$ choice. The compatibility with basepoint actions in the cited
theorem is what allows the cobordism map to restrict to these reduced groups.

The decoration-independence argument of \cite[Theorem~3.11 and Section~4]{Ef-TQFT},
with $\w$ and $\z$ interchanged for this choice of reduction, gives equality
on the $\uu$-multiple of the reduced source. More explicitly, if $F$ and
$F'$ are the maps associated with two decorations, then for every reduced
input $x$,
\[
\uu\bigl(F(x)-F'(x)\bigr)=F(\uu x)-F'(\uu x)=0.
\]
Their difference is therefore $\uu$-torsion in the target and vanishes
after passage to its torsion-free quotient. Moreover, bilinearity sends
an input with a torsion factor to torsion, so the quotient map is
well-defined. Thus the induced map
\begin{align*}
	\overline F_{\heeg,\spinct}:\Ical(Y_0,K_0,\spinc_0)\otimes_{\Abb}
	\Ical(Y_1,K_1,\spinc_1)\to \Ical(Y_2,K_2,\spinc_2),
\end{align*}
is well-defined and only depends on $\spinct$ and the surface $S_\heeg$ determined by $\heeg$ (and not its decoration).

When $\heeg$ is a grid triple, $Y_0,Y_1,Y_2\simeq S^3$ and there is only one $\SpinC$ structure on each boundary component. Moreover, each $\Ical(K_i)=\Ical(Y_i,K_i,\spinc_i)$ can be represented by a monomial ideal. This gives the maps
\begin{align*}
	\overline F_{\heeg,\spinct}:\Ical(K_0)\otimes_{\Abb}
	\Ical(K_1)\to \Ical(K_2).
\end{align*}

The absolute grading formula for the grid diagrams uses the formulas
\[
M_q(\x)=J(\x,\x)-2J(\x,q)+J(q,q)+1,\qquad q=\w,\z.
\]
For a grid of size \(n\) these are related to the symmetric cobordism
normalization of \cite{Zemke-gradings-revision} by
\[
\grw(\x)=M_{\w}(\x)+\frac{n-1}{2},\qquad
\grz(\x)=M_{\z}(\x)+\frac{n-1}{2}.
\]
We keep the conventional grid Alexander grading
\[A_{\rm grid}=\frac{M_{\w}-M_{\z}-n+1}2\] 
in Section~\ref{subsec:affine-alexander};
the symmetric multi-based Alexander grading is
\[\frac{\grw-\grz}2=A_{\rm grid}+\frac{n-1}2.\]

For a term from an index-zero triangle \(\psi\), its raw grid bidegree is
\begin{equation}\label{eq:raw-triangle-degree}
\begin{aligned}
\Delta_{\w}^{\rm grid}&=M_{\w}(\x_2)-M_{\w}(\x_0)-M_{\w}(\x_1)-2n_{\w}(\psi),\\
\Delta_{\z}^{\rm grid}&=M_{\z}(\x_2)-M_{\z}(\x_0)-M_{\z}(\x_1)-2n_{\z}(\psi).
\end{aligned}
\end{equation}
The cobordism degree, as computed by Zemke, is obtained by subtracting \((n-1)/2\) in each
coordinate, since there are two inputs and one output.
Write \(c_{\spinct}=\langle c_1(\spinct),H\rangle\), an odd integer.
Since \(H^2=-1\), we have \(c_1(\spinct)^2=-c_{\spinct}^2\).
For the midpoint decoration each colored subsurface is a union of \(n\)
disks with collars. Thus
\[
\chi(S_{\w})=\chi(S_{\z})=n,\qquad
\widetilde\chi(S_{\w})=\widetilde\chi(S_{\z})=n-\frac{3n}{2}=-\frac n2.
\]
Together with \(\chi(X)=0\), \(\sigma(X)=-1\), and
\([S_{\heeg}]=0\), the grading-change theorem of Zemke
\cite[Theorem~1.4]{Zemke-gradings-revision} gives
\begin{equation}\label{eq:bigrading-shift}
\deg f_{\heeg,\spinct}
=(d_{\spinct},d_{\spinct}),\qquad\text{where}\ \ \ 
d_{\spinct}=\frac{3-2n-c_{\spinct}^{2}}4.
\end{equation}
The same shift applies to a nonzero \(\z\)-blocked summand with these
absolute gradings. If \(\spinct=\spinct_0+mh\), 
\[
c_{\spinct}=c_{\spinct_0}+2m,\qquad
d_{\spinct}=d_{\spinct_0}-m c_{\spinct_0}-m^2.
\]
For the reduced map, the bidegree is \((d_{\spinct}+\rho_{\w},d_{\spinct}+\rho_{\z})\);
in the tensor choices above \(\rho_{\w}=-\rho_{\z}=\tfrac{n-1}2\).
With these choices, we obtain the map 
\begin{align*}
	&\overline F_{\heeg,\spinct}:
	\Ical(K_0)\otimes_{\Abb}\Ical(K_1)\longrightarrow\Ical(K_2),\qquad
	\deg(\overline F_{\heeg,\spinct})
	=\left(\frac{1-c_{\spinct}^2}{4},\frac{1-c_{\spinct}^2}{4}-n+1\right).
\end{align*}	

With the above degree normalization, we also have the monomial-ideal identifications
\[\Ical(K)\cong \Jcal(K)[\eta_w(K),\eta_z(K)],\qquad\text{for some}\ \ 
\eta(K)=(\eta_w(K),\eta_z(K))\in \Z^2, \] 
where \(\Jcal(K)\subset\Abb\) is a bigraded ideal with bigrading induced from $\Abb$ and \(\Jcal(K)[i,j]\) means that the degree of a monomial is increased by the bidegree \((i,j)\).

\begin{thm}\label{thm:ideal-inclusion}
Let \(K_0,K_1,K_2\) be the boundary knots of a centrally marked appropriate
triple grid. The bilinear homology map induces the homomorphisms of monomial ideals
\begin{align}
\overline F_{\heeg,\spinct}:\Ical(K_0)\otimes_{\Abb}\Ical(K_1)\to\Ical(K_2)
\qquad\text{and}\qquad
G_{\heeg,\spinct}:\Jcal(K_0)\otimes_{\Abb}\Jcal(K_1)\to\Jcal(K_2).
\end{align}
If \(\overline F_{\heeg,\spinct}\ne0\), then after
localization, both maps are multiplication by \(\uu^a\vv^b\), where 
\begin{align*}
&2a=-\eta_w(K_0)-\eta_w(K_1)+\eta_w(K_2)+\frac{c_{\spinct}^2-1}{4}
\qquad\text{and}\\
&2b=-\eta_z(K_0)-\eta_z(K_1)+\eta_z(K_2)+\frac{c_{\spinct}^2-1}{4}+n-1.	
\end{align*}	
In particular, in the common fraction field
\[
\uu^a\vv^b\Jcal(K_0)\Jcal(K_1)\subset\Jcal(K_2).
\]
\end{thm}
\begin{proof}
We have already proved most of the claims, and only the last part needs a justification.	
After tensoring with \(\mathbb F(u,v)\), each module is one-dimensional.
A nonzero map is multiplication by an element of \(\mathbb F(u,v)\). Homogeneity for the two independent variable degrees makes that element a Laurent monomial (with scalar coefficient one over \(\mathbb F_2\)). Let us assume that the monomial is $\uu^a\vv^b$. It then follows that $G_{\heeg,\spinct}=\uu^{a}\vv^{b}$ as a bigraded map. 
Accounting for the three ideal shifts shows that the bidegree of $G_{\heeg,\spinct}$ is
\begin{align*}
	\deg(G_{\heeg,\spinct})=\deg(\overline F_{\heeg,\spinct})+\left(\eta(K_0)+\eta(K_1)-\eta(K_2)\right)=(-2a,-2b).
\end{align*}	
This completes the proof.
\end{proof}
Negative exponents are permitted; the containment can equivalently be
written after clearing denominators. No containment is asserted for a
zero reduced map.

\subsection{The distinguished corner generators}\label{subsec:distinguished-generators}
Every ordinary grid diagram has canonical corner states determined by its
marked squares. Let \(\x^{NE}\), \(\x^{SW}\), \(\x^{NW}\), and \(\x^{SE}\)
denote the states occupying the northeast, southwest, northwest and southeast corners of all \(z\)-marked squares, respectively.

In the \(z\)-blocked complex, \(\x^{NE}\) and \(\x^{SW}\) are cycles;
see \cite[Proposition~6.3.1 and Remark~6.3.2]{OSS-book}.
When the distinguished class realizes the extremal non-torsion Alexander
grading in the appropriate knot complex, that grading determines \(\tau\).

To identify corners in all three ordered grids, use the orientation-preserving
coordinates
\begin{align*}
&(x,y)\ \ \text{for}\ \ (\alpha,\beta),\qquad
&&(\tfrac12-x-y,x)\ \ \text{for}\ \ (\beta,\gamma)
&&\text{and}
&&(x+y-\tfrac12,y)\ \ \text{for}\ \ (\alpha,\gamma).
\end{align*}
The first curve family is horizontal and the second vertical in each chart.
For a marking \(z\) in the \((i,j)\)-square, the triangle
\[
\mathcal T_{z}=U_{i,j}+U_{i-1,j}+U_{i,j-1}+H_{i,j}
\]
has vertices
\[
A_{i,j},\qquad B_{i,j-1}\qquad\text{and}\qquad C_{i-1,j}.
\]
They are respectively the northeast, northeast, and northwest corners of
the marked cells in these three charts. Its leg length is \(3/2\).
The triangles \(\mathcal T_{z}\) give a class
\[
\psi_{\z}\in\pi_2^0(\x_0^{NE},\x_1^{NE},\x_2^{NW}),\qquad
\widetilde\delta(\psi_{\z})=\frac32,\qquad
n_{\z}(\psi_{\z})=n,\quad n_{\w}(\psi_{\z})=0.
\]
The sheets are distinct translates with equal signed lengths, so Proposition~\ref{prop:sheet-separation} and Lemma~\ref{lem:nice-count} show that $\psi_{\z}$ is a nice contributing class. Put \(\spinct_{\z}=\spinc(\psi_{\z})\).

\begin{thm}\label{thm:top-generator-triangle}
For the unblocked, unreduced triangle map,
\[
\left\langle f^-_{\heeg,\spinct_{\z}}
(\x_0^{NE}\otimes\x_1^{NE}),\x_2^{NW}\right\rangle=\vv^n.
\]
There is exactly one holomorphically contributing class in this summand
with these three corners, namely \(\psi_{\z}\).
Its raw grid bidegree is \((-2,-2)\), and its cobordism bidegree is
\[
\left(-\frac{n+3}{2},-\frac{n+3}{2}\right).
\]
\end{thm}

\begin{proof}
The preceding construction and Lemma~\ref{lem:nice-count} give the
contribution \(\vv^n\).
If \(\psi\) is another contributing index-zero class with the same corners
and Spin$^c$ structure, homogeneity in \eqref{eq:raw-triangle-degree} gives
\(n_{\w}(\psi)=0\) and \(n_{\z}(\psi)=n\). The \(w\)-blocked version of the
embedding argument in Section~\ref{holo} makes all its sheets embedded,
with leg lengths less than \(n\).

At each \(A_{i,j}\) in the first corner state, the other two vertices
are fixed by its unique vertical and horizontal boundary curves:
they must be \(B_{i,j-1}\) and \(C_{i-1,j}\).
Let \(t\) denote the signed displacement from this \(A\)-corner to its
\(B\)-corner in the vertical direction. Then \(t\equiv-3/2\pmod n\).
Since \(|t|<n\), the only choices are \(-3/2\) and \(n-3/2\).
Their difference numbers are \(-3\) and \(2n-3\), respectively: a triangle
of signed leg \(t\) has \(\delta=2t\), as follows by counting its upper
and lower elementary triangles. If \(k\) sheets use the second choice,
\[
\delta(\psi)=-3n+2nk.
\]
Theorem~\ref{thm:delta-on-triangles} and the fixed Spin$^c$ structure force
\(k=0\). Thus every sheet is the specified \(\mathcal T_{z}\), proving
uniqueness of the contributing class and the coefficient formula.

The corner \(J\)-formulas give \(M_{\z}(\x^{NE})=1-n\) and
\(M_{\z}(\x^{NW})=0\). Hence its raw \(z\)-degree is
\(0-2(1-n)-2n=-2\). Equation~\eqref{eq:bigrading-shift}, using the
nullhomology of the surface, makes the \(w\)-degree equal to it.
The normalization change subtracts \((n-1)/2\) from each raw degree.
\end{proof}

\begin{cor}\label{cor:absolute-chern-coordinate}
For every triangle class in the centrally marked grid,
\[
\langle c_1(\spinc(\psi)),H\rangle=-\frac{\delta(\psi)}n.
\]
In particular \(c_{\spinct_0}=1\), \(c_{\spinct_{\z}}=3\), and
the unreduced cobordism degree of a summand represented by \(\psi\) is
\[
\left(\frac{3-2n-(\delta(\psi)/n)^2}{4},
      \frac{3-2n-(\delta(\psi)/n)^2}{4}\right).
\]
\end{cor}
\begin{proof}
In the charts above, the small-upper reference \(\psi_0\) has corners
\(\x_0^{NE},\x_1^{NW},\x_2^{NE}\) and \(n_{\z}=0\).
Its raw \(z\)-degree is zero. Thus its cobordism degree is
\(-(n-1)/2\), and \eqref{eq:bigrading-shift} gives \(c_0^2=1\).
The theorem gives \(c_{\spinct_{\z}}^2=9\).
Since \(\delta(\psi_{\z})-\delta(\psi_0)=-2n\),
Theorem~\ref{thm:delta-on-triangles} gives \(c_{\spinct_{\z}}=c_0+2\).
These equations force \(c_0=1\) and \(c_{\spinct_{\z}}=3\).
The general formula follows from the affine difference formula.
\end{proof}

\begin{rem}
The theorem is a statement about a chain coefficient. It does not assert
that its input states are unblocked cycles or that any selected basepoint
reduction is nonzero on torsion-free homology. In particular, no bound on
the least power of \(\uu\) occurring at \(\x^{NW}\) in a non-torsion
\(z\)-blocked cycle follows from this \(\vv^n\)-term alone.
The nonvanishing hypothesis of Theorem~\ref{thm:ideal-inclusion} remains
necessary.
\end{rem}

\section{A family of examples: affine grid triples}

\subsection{Affine grid knots and links}\label{subsec:affine-grid}

We begin with a family of ordinary grids which will be used throughout the examples.  

\begin{defi}\label{def:affine-knot}
	Fix $n\geq2$, let $p\in\mathbb Z/n\mathbb Z$ be a unit, and let $s\in\mathbb Z/n\mathbb Z$ be nonzero.  The \emph{affine grid} $G(n;p,s)$ is the ordinary $n\times n$ grid whose $z$- and $w$-markings in column $i$ lie in rows
	\[
	x_i\equiv n-1-pi\pmod n,
	\qquad
	y_i\equiv x_i-s\pmod n.
	\]
	We write $K(n;p,s)$ for the oriented link represented by this grid with the convention that $z=X$, $w=O$, vertical arcs run from $z$ to $w$, horizontal arcs from $w$ to $z$, and vertical arcs pass over horizontal arcs.
\end{defi}

All parameters following the semicolon are understood modulo $n$ unless an integer representative is explicitly specified.

\begin{lem}
	\label{lem:affine-components}
	The link $K(n;p,s)$ has
	\(\gcd(n,s)\)
	components.  In particular, $K(n;p,s)$ is a knot exactly when $\gcd(n,s)=1$.
\end{lem}

\begin{proof}
	Let $b=[p^{-1}s]_n\in\{1,\ldots,n-1\}$.	
	Since
	\[
	y_i\equiv x_i-s\equiv x_{i+b}\pmod n,
	\]
	following a vertical arc from $z_i$ to $w_i$ and then the horizontal arc returns to the $z$-marking in column $i+b$.  Thus the components are the cycles of translation by $b$ on $\mathbb Z/n\mathbb Z$.  Since multiplication by $p$ is invertible modulo $n$, $\gcd(n,b)=\gcd(n,s)$.
\end{proof}

The elementary grid symmetries will be useful below.  Interchanging the $z$- and $w$-markings reverses the orientation of every component and replaces $s$ by $-s$.  Reflection in a horizontal axis takes the mirror and replaces $(p,s)$ by $(-p,-s)$.  Thus
\[
K(n;p,-s)=-K(n;p,s),
\qquad
K(n;-p,-s)=\overline{K(n;p,s)}
\]
up to the harmless cyclic row and column translations needed to restore the chosen constant term $n-1$.

Two further elementary symmetries make the ordered grid charts below
explicit. Rotation through \(180^\circ\), followed by translations, takes
\(G(n;p,s)\) to \(G(n;p,-s)\) without changing its oriented link type.
Thus these affine links are isotopic to their simultaneous orientation
reversals. Transposing a planar grid exchanges the crossing convention
as well as the coordinate axes; the combined operation preserves the
ambient orientation and reverses all link orientations. Solving for its
row permutation gives
\[
K(n;p,s)\cong K(n;p^{-1},p^{-1}s).
\]
Here simultaneous link reversal has been absorbed using the preceding
rotation symmetry. These identities concern link types, not equality of
particular corner states.

\begin{prop}
	\label{prop:affine-torus}
	Suppose that $p,q$ are positive integers, \(s=pq<n\) and $\gcd(p,n)=1$.
	Then
	\[K(n;p,pq)=T_{q,n-pq}.\]
\end{prop}

\begin{proof}
Put \(m=n-pq>0\). On the toroidal grid choose each vertical arc from
\(z_i\) to \(w_i\) to run upwards a distance \(m\), and each horizontal
arc from \(w_i\) to \(z_{i+q}\) to run rightwards a distance \(q\).
These are valid bridge representatives: changing to the complementary
arc in a row or column uses its compressing disk and does not change
the grid link. We check that their union on the torus is embedded.

A possible interior intersection of the horizontal arc beginning in
column \(i\) with a vertical arc must occur in column \(i+k\), where
\(1\le k<q\). The height of that horizontal arc relative to the bottom
\(z_{i+k}\) of the vertical arc is, modulo \(n\),
\[
m+pk=n-p(q-k).
\]
This lies strictly between \(m\) and \(n\), so it is outside that vertical
arc. Parallel arcs also have disjoint interiors. Consequently the entire
staircase link is embedded on the Heegaard torus, and the usual small
pushes of its vertical and horizontal arcs can be undone without crossings.

Its total horizontal and vertical winding numbers are \(q\) and \(m\):
the \(n\) horizontal arcs have total length \(nq\), and the vertical arcs
have total length \(nm\). Its components are disjoint essential curves
and therefore parallel. If \(g=\gcd(q,m)=\gcd(n,pq)\), each has primitive
slope \((q/g,m/g)\). Thus the link is the standard \(T_{q,m}\), with the
orientation prescribed in Definition~\ref{def:affine-knot}.
\end{proof}

\begin{rem}\label{rem:affine-torus-recognition}
	Proposition~\ref{prop:affine-torus} is deliberately an integer statement, rather than merely a congruence claim.  A congruence $pq\equiv s\pmod n$ does not in general identify $K(n;p,s)$ with $T_{q,n-s}$; wrapping in the planar grid may change the knot type.  The symmetries above nevertheless give many further torus identifications by first replacing $(p,s)$ by one of the equivalent parameter pairs and then applying Proposition~\ref{prop:affine-torus}.
\end{rem}

\subsection{The distinguished bigradings of affine knots}\label{subsec:affine-alexander}

The affine form makes an otherwise separate calculation completely explicit: the bigrading of the distinguished northeast-$z$ state.  This calculation remains useful independently of the triangle-map discussion.  

Assume $1\le p<n$, $\gcd(p,n)=1$, and $0<s<n$.
Let $\x^+_{n;p,s}$ denote the state occupying the northeast corner of every $z$-marked square of $G(n;p,s)$.  As before, set \(b=[p^{-1}s]_n\in\{1,\ldots,n-1\}\)  and define

\begin{align*}
	&\Lambda(n;p,s)
	=
	\frac{(b-1)(pb-s)}{n}
	-2\sum_{j=1}^{b-1}\left\lfloor\frac{pj}{n}\right\rfloor\qquad\text{and}\\
	&\mu(n;p,s)=-(b+1)(n-s-1)+n\Lambda(n;p,s).
\end{align*}
We may then extend the definition and set \[\Lambda(n;p,s)=\Lambda(n;[p]_n,[s]_n)\qquad\text{and}\qquad
\mu(n;p,s)=\mu(n;[p]_n,[s]_n)\]
for arbitrary integers $p$ and $s$, where $\gcd(n,p)=1$ and $s$ is not a multiple of $n$.
The apparently asymmetric formula is invariant under cyclic row and column translations
of the affine grid.

\begin{thm}
	\label{thm:affine-distinguished-bigrading}
	For the northeast-$z$ state of $G(n;p,s)$ one has
	\[
	\bigl(M_{\w}(\x^+_{n;p,s}),
	M_{\z}(\x^+_{n;p,s})\bigr)
	=\bigl(\mu(n;p,s),1-n\bigr).
	\]
	Equivalently, the grid Alexander grading $A_{\rm grid}$ of the aforementioned grid state is 
	given by
	\[
	A_{\rm grid}(\x^+_{n;p,s})=\frac{\mu(n;p,s)}2.
	\]
\end{thm}

\begin{proof}
	The universal northeast-corner calculation gives
	\[
	M_{\z}(\x^+)=1-n.
	\]
	For the first Maslov grading, evaluate the standard $J$-formula
	\[
	M_{\w}(\x)=J(\x,\x)-2J(\x,\w)
	+J(\w,\w)+1
	\]
	on the affine permutations. 
	
For clarity, put \(r_i=[n-1-pi]_n\) and \(t_i=[r_i-s]_n\), using
\(0\le i<n\). In the fundamental square \((0,n]^2\), the northeast
state has coordinates \((i+1,r_i+1)\). Expanding the three \(J\)-terms gives
\[
\begin{split}
M_{\w}(\x^+)=1-\sum_i\mathbf1_{\{t_i\le r_i\}}
+\sum_{i<j}\bigl(&\mathbf1_{\{r_i<r_j\}}+\mathbf1_{\{t_i<t_j\}}
-\mathbf1_{\{r_i<t_j\}}-\mathbf1_{\{t_i\le r_j\}}\bigr).
\end{split}
\]
The single sum is \(n-s\). For the pair sum, use \(t_i=r_{i+b}\),
with indices read cyclically, where \(b=[p^{-1}s]_n\).
Cancel the pairs common to the two cyclic orders. The remaining pairs
give
\[
\sum_{i<j}\bigl(\mathbf1_{\{r_i<r_j\}}+\mathbf1_{\{t_i<t_j\}}
-\mathbf1_{\{r_i<t_j\}}-\mathbf1_{\{t_i\le r_j\}}\bigr)
=-b(n-s-1)+n\sum_{j=1}^{b-1}\mathbf1_{\{[pj]_n>s\}}.
\]
Thus the explicit \(J\)-expansion yields

	\[M_{\w}(\x^+)=
	-(b+1)(n-s-1)+n\#\{1\le j<b:[pj]_n>s\}.\]
	It remains only to rewrite the last count.  Write $pb=s+n\kappa$, for a non-negative integer $\kappa$. Then, for
	$1\le j<b$ one has
	\begin{align*}
		&\mathbf 1_{\{[pj]_n>s\}}
		=
		\left\lfloor\frac{pj+n-1-s}{n}\right\rfloor
		-
		\left\lfloor\frac{pj}{n}\right\rfloor\qquad
		\Rightarrow\\
		&\#\{1\le j<b:[pj]_n>s\}
		=\sum_{j=1}^{b-1}\left\lfloor\frac{p(b-j)+n-1-s}{n}\right\rfloor-
		\sum_{j=1}^{b-1}\left\lfloor\frac{pj}{n}\right\rfloor
		=(b-1)\kappa
		-2\sum_{j=1}^{b-1}\left\lfloor\frac{pj}{n}\right\rfloor,
	\end{align*}
	which is $\Lambda(n;p,s)$.  Finally, the grid identity
	\[
	A_{\rm grid}(\x)=\frac{M_{\w}(\x)-M_{\z}(\x)}2-\frac{n-1}{2}
	\]
	reduces here to $A_{\rm grid}(\x^+)=M_{\w}(\x^+)/2$.
\end{proof}

\subsection{Affine triple grids}\label{subsec:affine-cobordism}

We now return to the construction of triple grids.  Start with the markings of $G(n;p,s)$ in the $\alpha\beta$ diagram and assume in addition that $p-1$ is a unit modulo $n$.  Put
\[
r=(p-1)^*\in(\mathbb Z/n\mathbb Z)^\times.
\]
The order of the two attaching-curve families matters when the other two pairwise diagrams are read as ordinary grids: for an ordered pair $(\lambda,\mu)$ we draw the first family horizontally, the second vertically, and let the second family pass over the first.

\begin{thm}
	\label{thm:affine-cobordism}
	With the convention above, the three sub-diagrams of $\heeg$ and the three corresponding boundary links of the pointed grid triple are
	\begin{align*}
		&\heeg_0:\ K(n;p,s),
		&&\heeg_1:\ K(n;1-p,-s)
		&&\text{and}
		&&\heeg_2:\ K(n;-pr,-sr).
	\end{align*}
	Consequently the triple grid determines the cobordism maps
	\begin{align*}
		&K(n;p,s)\coprod K(n;1-p,-s)\longrightarrow K(n;-pr,-sr)\qquad\text{and}\\
		&K(n;pr,sr)\coprod K(n;p,s)\longrightarrow K(n;p-1,s)
	\end{align*}
\end{thm}

\begin{proof}
Use the three orientation-preserving charts in
Section~\ref{subsec:distinguished-generators}. In the \(\beta\gamma\) chart
\((X,Y)=(\tfrac12-x-y,x)\), the centers of the \(z\)-marked cells satisfy
\[
X\equiv(p-1)i+\tfrac12\pmod n,\qquad Y=i+\tfrac12.
\]
Solving for the row as a function of the new column gives the parameters
\((-r,rs)\), up to a cyclic row translation. By the transpose symmetry,
\[
K(n;-r,rs)\cong K(n;1-p,-s).
\]
In the \(\alpha\gamma\) chart
\((X,Y)=(x+y-\tfrac12,y)\), the new column is
\((1-p)i-1\) modulo \(n\). Solving the row equation gives slope
parameter \(-pr\). The \(w\)-row at a fixed new column differs from the
\(z\)-row by \(rs\), so its separation parameter is \(-sr\).
This proves the three link identifications.

The first cobordism is the ordered triangle map
\[C(\alpha,\beta)\otimes C(\beta,\gamma)\to C(\alpha,\gamma).\]
Cyclically ordering the same triple as \((\gamma,\alpha,\beta)\) gives
\[C(\gamma,\alpha)\otimes C(\alpha,\beta)\to C(\gamma,\beta).\]
Reversing an ordered pair reverses the ambient three-manifold orientation,
so its link is mirrored (any simultaneous component reversal is harmless
for this affine family). The mirror identities above yield the second
displayed cobordism. Both are in \(\CPbar\) with three balls removed,
with their indicated standard incoming and outgoing boundaries.
\end{proof}

\begin{cor}
	\label{cor:p2-family}
	Let $n$ be odd and $q\ge1$, let $s=2q<n$, and put $m=n-2q$.  Then $G(n;2,2q)$ corresponds to the decorated cobordisms
	\[
	T_{q,m}\coprod T_{2q,-m}\longrightarrow T_{q,-m}\qquad\text{and}\qquad
	T_{q,m}\coprod T_{q,m}\longrightarrow T_{2q,m}.
	\]
	Both are hosted by $\CPbar$ minus three $4$-balls.
	When $\gcd(q,m)=1$, all three boundary links are torus knots.
\end{cor}

\begin{proof}
	Here $r=1$ modulo $n$.  Proposition~\ref{prop:affine-torus} gives
	$K(n;2,2q)=T_{q,m}$.  Moreover
	\begin{align*}
		\overline{K(n;-1,-2q)}=K(n;1,2q)=T_{2q,m}
		\qquad\text{and}\qquad
		\overline{K(n;-2,-2q)}=K(n;2,2q)=T_{q,m}.
	\end{align*}
	Now apply Theorem~\ref{thm:affine-cobordism}.
\end{proof}

A second useful family shows that $p=2$ is not the only case in which all three boundary knots can be identified explicitly.

\begin{cor}
	\label{cor:p3-family}
	Assume $n>6$ and $\gcd(n,6)=1$.  Then the choice $(p,s)=(3,6)$ gives
	\[
	T_{2,n-6}\coprod T_{3,6-n}
	\longrightarrow T_{2,3}\qquad\text{and}\qquad
	T_{2,n-6}\coprod T_{2,-3}
	\longrightarrow T_{3,n-6}.
	\]
\end{cor}

\begin{proof}
	The first knot is $K(n;3,6)=T_{2,n-6}$.  The second is
	\[
	K(n;-2,-6)=\overline{K(n;2,6)}=T_{3,6-n}.
	\]
	Since $r=(2)^*=(1-n)/2$ modulo $n$, we may write
	\[
	pr=\frac{3-n}{2},
	\qquad
	sr=3-n.
	\]
	The output of the native cobordism has parameters
	\[
	-pr=\frac{n-3}{2},
	\qquad
	-sr=n-3,
	\]
	and $n-3=2((n-3)/2)<n$.  Proposition~\ref{prop:affine-torus} therefore identifies the output with $T_{2,3}$.
\end{proof}

For the affine family put \(g=\gcd(n,s)\). Each boundary link has \(g\)
components. The two-step moves between \(\z\)-markings in the bridge graph
translate the column by \(p^{-1}s\) or \(rs\). Their orbits are precisely
the residue classes modulo \(g\), so the surface has \(g\) components.
Each has three boundary circles and \(2n/g\) disks and \(3n/g\) bands,
hence genus \((n/g-1)/2\). In the knot cases of
Corollaries~\ref{cor:p2-family} and \ref{cor:p3-family}, the surface is
connected of genus \((n-1)/2\). All these surfaces are relatively
nullhomologous by Proposition~\ref{prop:nullhomologous-grid-surface}.

\section{Tool Use}
Assistance from GPT-6 has been crucial in certain parts of the paper. Most importantly, the current proof of Theorem~\ref{thm:holomorphic-domains-are-nice} was suggested to replace a much longer proof by the author. The separation criteria of Section~\ref{subsec:sheet-separation} also grew from discussions with GPT-6.  It also helped in editing the paper for better flow and proofreading, which resulted in finding several earlier mistakes in bigrading computations.

\end{document}

%% file: Tri-Grid.pdf_tex
%% Creator: Inkscape 1.2.1 (9c6d41e, 2022-07-14), www.inkscape.org
%% Accompanies image file 'Tri-Grid.pdf' (pdf, eps, ps)
\begingroup%
  \makeatletter%
  \providecommand\color[2][]{%
    \errmessage{(Inkscape) Color is used for the text in Inkscape, but the package 'color.sty' is not loaded}%
    \renewcommand\color[2][]{}%
  }%
  \providecommand\transparent[1]{%
    \errmessage{(Inkscape) Transparency is used (non-zero) for the text in Inkscape, but the package 'transparent.sty' is not loaded}%
    \renewcommand\transparent[1]{}%
  }%
  \providecommand\rotatebox[2]{#2}%
  \newcommand*\fsize{\dimexpr\f@size pt\relax}%
  \newcommand*\lineheight[1]{\fontsize{\fsize}{#1\fsize}\selectfont}%
  \ifx\svgwidth\undefined%
    \setlength{\unitlength}{359.25000208bp}%
    \ifx\svgscale\undefined%
      \relax%
    \else%
      \setlength{\unitlength}{\unitlength * \real{\svgscale}}%
    \fi%
  \else%
    \setlength{\unitlength}{\svgwidth}%
  \fi%
  \global\let\svgwidth\undefined%
  \global\let\svgscale\undefined%
  \makeatother%
  \begin{picture}(1,0.99965964)%
    \lineheight{1}%
    \setlength\tabcolsep{0pt}%
    \put(0,0){\includegraphics[width=\unitlength,page=1]{Tri-Grid.pdf}}%
  \end{picture}%
\endgroup%

%% file: G0G1G2.pdf_tex
%% Creator: Inkscape 1.2.1 (9c6d41e, 2022-07-14), www.inkscape.org
%% Accompanies image file 'G0G1G2.pdf' (pdf, eps, ps)
\begingroup%
  \makeatletter%
  \providecommand\color[2][]{%
    \errmessage{(Inkscape) Color is used for the text in Inkscape, but the package 'color.sty' is not loaded}%
    \renewcommand\color[2][]{}%
  }%
  \providecommand\transparent[1]{%
    \errmessage{(Inkscape) Transparency is used (non-zero) for the text in Inkscape, but the package 'transparent.sty' is not loaded}%
    \renewcommand\transparent[1]{}%
  }%
  \providecommand\rotatebox[2]{#2}%
  \newcommand*\fsize{\dimexpr\f@size pt\relax}%
  \newcommand*\lineheight[1]{\fontsize{\fsize}{#1\fsize}\selectfont}%
  \ifx\svgwidth\undefined%
    \setlength{\unitlength}{1025.12132034bp}%
    \ifx\svgscale\undefined%
      \relax%
    \else%
      \setlength{\unitlength}{\unitlength * \real{\svgscale}}%
    \fi%
  \else%
    \setlength{\unitlength}{\svgwidth}%
  \fi%
  \global\let\svgwidth\undefined%
  \global\let\svgscale\undefined%
  \makeatother%
  \begin{picture}(1,0.90424914)%
    \lineheight{1}%
    \setlength\tabcolsep{0pt}%
    \put(0,0){\includegraphics[width=\unitlength,page=1]{G0G1G2.pdf}}%
    \put(0.14778739,0.44711172){\color[rgb]{0.00392157,0.00392157,0.01176471}\makebox(0,0)[lt]{\lineheight{1.25}\smash{\begin{tabular}[t]{l}(a)\end{tabular}}}}%
    \put(0.14778739,0.00813929){\color[rgb]{0.00392157,0.00392157,0.01176471}\makebox(0,0)[lt]{\lineheight{1.25}\smash{\begin{tabular}[t]{l}(b)\end{tabular}}}}%
    \put(0.5721274,0.44711172){\color[rgb]{0.00392157,0.00392157,0.01176471}\makebox(0,0)[lt]{\lineheight{1.25}\smash{\begin{tabular}[t]{l}(c)\end{tabular}}}}%
    \put(0.5721274,0.00813929){\color[rgb]{0.00392157,0.00392157,0.01176471}\makebox(0,0)[lt]{\lineheight{1.25}\smash{\begin{tabular}[t]{l}(d)\end{tabular}}}}%
  \end{picture}%
\endgroup%

%% file: Diff-Number.pdf_tex
%% Creator: Inkscape 1.2.1 (9c6d41e, 2022-07-14), www.inkscape.org
%% Accompanies image file 'Diff-Number.pdf' (pdf, eps, ps)
\begingroup%
  \makeatletter%
  \providecommand\color[2][]{%
    \errmessage{(Inkscape) Color is used for the text in Inkscape, but the package 'color.sty' is not loaded}%
    \renewcommand\color[2][]{}%
  }%
  \providecommand\transparent[1]{%
    \errmessage{(Inkscape) Transparency is used (non-zero) for the text in Inkscape, but the package 'transparent.sty' is not loaded}%
    \renewcommand\transparent[1]{}%
  }%
  \providecommand\rotatebox[2]{#2}%
  \newcommand*\fsize{\dimexpr\f@size pt\relax}%
  \newcommand*\lineheight[1]{\fontsize{\fsize}{#1\fsize}\selectfont}%
  \ifx\svgwidth\undefined%
    \setlength{\unitlength}{377.25000941bp}%
    \ifx\svgscale\undefined%
      \relax%
    \else%
      \setlength{\unitlength}{\unitlength * \real{\svgscale}}%
    \fi%
  \else%
    \setlength{\unitlength}{\svgwidth}%
  \fi%
  \global\let\svgwidth\undefined%
  \global\let\svgscale\undefined%
  \makeatother%
  \begin{picture}(1,0.99602846)%
    \lineheight{1}%
    \setlength\tabcolsep{0pt}%
    \put(0,0){\includegraphics[width=\unitlength,page=1]{Diff-Number.pdf}}%
  \end{picture}%
\endgroup%